\documentclass[a4paper]{article}

\usepackage[bottom=2.8cm,top=2.8cm,left=2.8cm,right=2.8cm]{geometry}
\usepackage[OT1]{fontenc}
\usepackage{color,latexcolors}
\usepackage{soul}
\usepackage{amsthm,amsmath,graphicx,latexsym,amssymb,amscd,amsfonts,enumerate,bm}
\usepackage[colorlinks,citecolor=blue,linkcolor=red,urlcolor=magenta]{hyperref}
\usepackage{cancel,comment}
\usepackage[normalem]{ulem}

\usepackage{mathtools}
\usepackage{dsfont}

\newtheorem{theorem}{Theorem}[section]
\newtheorem{proposition}[theorem]{Proposition}
\newtheorem{lemma}[theorem]{Lemma}
\newtheorem{corollary}[theorem]{Corollary}

\theoremstyle{definition}
\newtheorem{remark}[theorem]{Remark}

\numberwithin{equation}{section}

\allowdisplaybreaks

\newcommand{\F}{\mathcal{F}}
\newcommand{\MM}{\mathcal{M}}
\def\E{\mathbb{E}}
\def\R{\mathbb{R}}
\def\P{\mathbb{P}}

\def\dd{\mathrm{d}}
\newcommand{\wt}{\widetilde}
\newcommand{\eps}{\varepsilon}
\newcommand{\N}{\mathbb{N}}
\newcommand{\1}{\mathds{1}}

\begin{document}

\title{\textbf{\large 
{Sharp estimates for incremental moments of inverse subordinators and approximation errors for time-changed stochastic differential equations}}}

\author{Ivan Bio\v{c}i\'c\\ 
{\small Department of Mathematics, Faculty of Science, University of Zagreb, Zagreb, Croatia.}\\ 
{\small Email: \texttt{ivan.biocic@math.hr}} 
\and 
Sixian Jin\\
{\small Department of Mathematics, California State University San Marcos, San Marcos, CA, USA.}\\
{\small Email: \texttt{sjin@csusm.edu}}
\and 
Kei Kobayashi\thanks{Corresponding author.}\\
{\small Department of Mathematics, Fordham University, New York, NY, USA.}\\	 {\small Email: \texttt{kkobayashi5@fordham.edu}}}

\maketitle                                                       
\vspace{-2mm}

\begin{abstract}
This paper establishes precise upper and lower bounds for the incremental moments of general inverse subordinators. Such processes are used as random time changes to model subdiffusive phenomena. The upper and lower bounds provide the foundation for a general framework for analyzing numerical approximation of time-changed stochastic differential equations on deterministic time grids. 
We derive an upper bound for the approximation error and further prove that the resulting order of convergence cannot be improved in general.
\vspace{2mm}

   \noindent

\noindent\textbf{Keywords:} inverse subordinator; time-changed stochastic differential equations; numerical approximation; strong error estimates; anomalous diffusion.\vspace{2mm}

\noindent
\textbf{MSC (2020)}: Primary 60H35; Secondary 65C30, 60H10, 60K50
\end{abstract}

\section{Introduction}
Let $B=(B_t)_{t\ge0}$ be a standard Brownian motion and $E=(E_t)_{t\ge0}$ be the inverse of a stable subordinator $D=(D_t)_{t\ge0}$ with index $\alpha\in(0,1)$, which is independent of $B$. The composition $B\circ E=(B_{E_t})_{t\ge0}$, which we call a \textit{time-changed Brownian motion}, and its various generalizations have been widely employed to model subdiffusions arising in many different areas of science.
The time-changed Brownian motion is non-Markovian,
and 
its variance function is $\E[B_{E_t}^2]=t^\alpha/\Gamma(1+\alpha)$,
which shows that on large time scales
particles represented by $B\circ E$ diffuse more slowly than the regular Brownian particles.  
 Moreover, the densities of $B\circ E$ satisfy the time-fractional heat equation $\partial_t^\alpha p(t,x) =\frac12\Delta p(t,x)$, with $\partial_t^\alpha$ denoting the Caputo fractional derivative of order $\alpha$. The fractional derivative $\partial_t^\alpha$ is a non-local operator, which indicates the non-Markovian nature of the time-changed Brownian motion.
Stochastic differential equations driven by $B\circ E$ and their associated fractional-order partial differential equations have been investigated.
Such stochastic differential equations are ideal models for capturing memory effects and subdiffusive behaviors observed in various complex real-world systems. For details of these accounts, see \cite{HKU-book} and references therein.
 
The main purpose of this paper is two-fold: (i) to establish precise upper and lower bounds for incremental moments of general inverse subordinators, and (ii) to apply these bounds to the error analysis of numerical approximation of solutions to stochastic differential equations driven by the time-changed Brownian motion $B\circ E$, with $E$ being a general inverse subordinator.

To be more precise, our first main result, Theorem \ref{Theorem_increments}, generalizes \cite[Lemma 2.1]{Gajda-Magdziarz_SPL} on inverse \textit{stable} subordinators to a wide class of general inverse subordinators. \cite[Lemma 2.1]{Gajda-Magdziarz_SPL} states that if $E$ is the inverse of an $\alpha$-stable subordinator $D$, then for any $n\in\N$ and $0\le a<b<\infty$,  
\begin{align}\label{Eincrements-stable}
	\Bigl(\frac{b-a}{b}\Bigr)^{1-\alpha}\frac{n!(b-a)^{n\alpha}}{\Gamma((n-1)\alpha+2)\Gamma(\alpha)}
	\le \E\left[ (E_b-E_a)^n\right]
	\le \Bigl(\frac{b-a}{b}\Bigr)^{1-\alpha}\frac{n!(b-a)^{n\alpha}}{\Gamma(n\alpha+1)}.
\end{align}
This allows us to infer the magnitudes of incremental moments of $E$ when equidistant deterministic time steps $0=t_0<t_1<t_2<\ldots$ with $\Delta\coloneqq t_{i+1}-t_i$ are taken. Indeed, \eqref{Eincrements-stable} implies that for small $\Delta$ and any $i\in\mathbb{N}$, 
\[
	\E[(E_{t_{i+1}}-E_{t_i})^n]\asymp t_{i+1}^{\alpha-1}\Delta^{1+(n-1)\alpha},
\]
where we use the notation $f\asymp g$ to indicate that there exists a positive constant $c$ such that $cf(x)\le g(x)\le f(x)/c$ for all $x$ in an appropriate domain.
In contrast, building upon the ideas presented in \cite{Lageras} and \cite{VeilletteTaqqu-letters}, Theorem \ref{Theorem_increments} establishes precise upper and lower bounds of incremental moments for a wide class of general inverse subordinators.  
In the most general setting of Theorem \ref{Theorem_increments}, we will express the bounds in terms of the renewal function $U(t)\coloneqq \E[E_t]$ (see \eqref{U-def}). While the non-Markovian, non-self-similar nature of general inverse subordinators requires careful, non-trivial analysis, the obtained expressions have a clear benefit; unlike \eqref{Eincrements-stable}, they explicitly show the inevitable connection of the incremental moments to potential theory. On the other hand, if the underlying Laplace exponent $\phi$ satisfies appropriate weak scaling conditions, the latter bounds reduce to forms similar to those in \eqref{Eincrements-stable}. We will present multiple forms of upper and lower bounds across various settings to help the reader utilize them in different contexts in which they are interested. 

Although various applications of Theorem \ref{Theorem_increments} are conceivable, this paper particularly focuses on applications to the numerical analysis of the following two classes of stochastic differential equations (SDEs) driven by the time-changed Brownian motion:
\begin{align}\label{403}
	dX_t=F(t,X_t)\,dE_t+	G(t,X_t)\,dB_{E_t},
\end{align}
and 
\begin{align}\label{402}
	dX_t=F(E_t,X_t)\,dE_t+G(E_t,X_t)\,dB_{E_t}.
\end{align}
Investigations into numerical approximations of the latter SDEs \eqref{403} and \eqref{402} with general Lipschitz coefficients, along with their orders of strong convergence, started within the past decade with papers \cite{JinKobayashi2021,JumKobayashi} and \cite{JinKobayashi2019}, respectively, which built upon the idea presented in \cite{Magdziarz_SPA}. Since then, similar investigations have been carried out for the time-changed SDEs and their generalizations, with different schemes and different assumptions imposed on the coefficients. See e.g., \cite{ChenLiuWu2026,DengLiu2020,Shen+2023} for convergence studies and \cite{NaneNi2017,NaneNi2018} for related stability analysis of time-changed SDEs. 
In these works concerning strong convergence, the approximation scheme involved first simulating the underlying subordinator $D$ at deterministic time points $0=t_0<t_1<t_2<\ldots$ and then discretizing the inverse subordinator $E$ at the random time points $\tau_i\coloneqq D_{t_i}$, $i\ge 0$. The advantage of using such a scheme is that it makes the increments of $E$ at the random time points non-random, in the sense that $E_{\tau_{i+1}}-E_{\tau_i}=t_{i+1}-t_i$, which results in a cleaner error analysis. On the other hand, the scheme may be considered unnatural since random (rather than deterministic) time points are employed to discretize the SDE. 
 In contrast, the recent paper \cite{BCT2026} introduced a scheme that allows the inverse subordinator $E$ and SDE \eqref{402} to be discretized at deterministic time points $0=t_0<t_1<t_2<\ldots$ directly, while still relying on a convergence result in \cite{JumKobayashi} to help obtain the order of strong convergence. 

This paper assumes that one already has an approximation scheme in hand to discretize the inverse subordinator $E$ at deterministic time points (as in the paper \cite{BCT2026}), while maintaining the spirit of the error analysis conducted in the papers \cite{JinKobayashi2019,JinKobayashi2021,JumKobayashi}. 
{Theorems \ref{thm:strong-approx} and \ref{thm:synchro-strong} establish the rates of the strong $L^2$ errors of our approximation schemes for SDEs \eqref{403} and \eqref{402}, respectively. 
The rate functions take the forms
$\Delta^{2\rho}+\phi(1/\Delta)^{-(1-\varepsilon)}$ 
 and 
 $\phi(1/\Delta)^{-(\min\{1,2\rho\}-\varepsilon)}$,
 respectively, where $\varepsilon>0$ is an arbitrarily small constant, $\phi$ denotes the Laplace exponent of the underlying subordinator, and $\rho\in(0,1]$ controls the regularity of the coefficients $F$ and $G$ with respect to the time variable. It is further verified in Theorem \ref{prop:sharp-piecewise} that the above rate functions cannot be generally improved to the forms with $\varepsilon=0$, which indicates the sharpness of our error estimates in Theorems \ref{thm:strong-approx} and \ref{thm:synchro-strong}.}

{To obtain the sharp error bounds,} 
precise estimates of the incremental moments of the inverse subordinator $E$ in Theorem \ref{Theorem_increments} are indispensable. 
In fact, the use of the well-known simple estimate $\E\left[ (E_b-E_a)^n\right]\le \E\left[ E_{b-a}^n\right]$ does not help achieve the rates of convergence established in Theorems \ref{thm:strong-approx} and \ref{thm:synchro-strong}.

{Another significant idea to be employed for precise analysis of the magnitude of the approximation error is \textit{exponential tilting} of the conditional error process, which 
substantially enlarges the class of inverse subordinators for which our error analysis is valid and also considerably improves a result established in \cite[Theorem 7]{JinKobayashi2019}; see Remark \ref{r:compare-deterministic-random} for details of this account. 
For classical It\^o SDEs (without time changes), exponential tilting is unnecessary, and 
to our knowledge, its use in this particular context of numerical approximation of time-changed SDEs is completely new.}

Let us emphasize that part of the purpose of this paper is to provide \textit{robust} machinery for the error analysis of approximation schemes for general time-changed SDEs when the time steps are naturally taken to be deterministic, rather than to identify the orders of convergence for different kinds of SDE coefficients and approximation schemes. To this end, we present our results and discussions in the simplest setting; throughout the paper, the coefficients of time-changed SDEs are assumed to be Lipschitz, and only Euler--Maruyama-type approximation schemes are considered. In addition, even though it is possible to allow the coefficients in SDEs \eqref{403} and \eqref{402} to take the merged forms $F(t,E_t,X_t)$ and $G(t,E_t,X_t)$, we avoid choosing this option since it would not only complicate the presentation but also obscure some aspects of the analysis of each SDE. For example, the duality principle established in \cite{Kobayashi2011} allows us to analyze SDE \eqref{402} by connecting it to the classical It\^o SDE without the random time change, whereas that is not the case for SDE \eqref{403}. 

Let us also note that a statement similar in spirit to Theorem \ref{thm:strong-approx}
has appeared in \cite[Theorem 3.1]{Zuo2026} under the assumption that the random time change is an inverse $\alpha$-\textit{stable} subordinator with $\alpha\in(\frac 12,1$). 
However, Theorem \ref{thm:strong-approx} is applicable to a much broader class of inverse subordinators, including the inverses of $\alpha$-stable subordinators with $\alpha\in(0,1)$ and general non-stable subordinators. Our approach also reveals a direct connection with the renewal function and also provides a remedy for the unjustified statement used in the proof of \cite[Theorem 3.1]{Zuo2026} in analyzing the error associated with the subdiffusion coefficient $G$. See Remark \ref{rem:Zuo} for details.   

The rest of the paper is organized as follows. Section \ref{s:increment_moment} defines inverse subordinators and provides results associated with their moments; in particular, Theorem \ref{Theorem_increments} establishes precise estimates of incremental moments.  Section \ref{s:approx-1} clarifies the meaning of the time-changed SDE \eqref{403} and its Euler--Maruyama-type approximation scheme, as well as derives the order of strong convergence of the approximation scheme in Theorem \ref{thm:strong-approx}. In contrast, Section \ref{s:approx-2} focuses on SDE \eqref{402} and establishes convergence results in Theorems \ref{thm:synchro-strong} and \ref{thm:approx-synchr-conti}. 
Section \ref{s:remarks} provides several noteworthy remarks about the statements derived in Sections \ref{s:approx-1} and \ref{s:approx-2}, including remarks that compare the convergence orders obtained in Theorems \ref{thm:strong-approx}, \ref{thm:synchro-strong}, and \ref{thm:approx-synchr-conti} with each other and with the existing statements in the literature.   
Lastly, Section \ref{s:sharp} provides an example showing that the convergence orders established in Theorems \ref{thm:strong-approx} and \ref{thm:synchro-strong} \textit{cannot} be generally improved, thereby confirming the sharpness of our error analysis.

\section{Moments associated with inverse subordinators}
\label{s:increment_moment}

This section establishes some important results on moments associated with general inverse subordinators. The first subsection derives upper and lower bounds for incremental moments of inverse subordinators (Theorem \ref{Theorem_increments}), which will play a crucial role in Sections \ref{s:approx-1} and \ref{s:approx-2} while analyzing the strong $L^2$ error of numerical approximation of time-changed SDEs of the forms \eqref{403} and \eqref{402}, respectively.  
On the other hand, the next subsection is devoted to establishing general criteria for the existence and non-existence of the exponential moments of powers of inverse subordinators (Theorem \ref{theorem_expo_moment}). 
Let us emphasize that the results to be presented in this section cover a broad class of inverse subordinators and substantially generalize the existing results available in the literature, including \cite[Lemma 2.1]{Gajda-Magdziarz_SPL} and \cite[Theorem 1]{JinKobayashi2019}. In particular, although Theorem \ref{theorem_expo_moment} is not utilized in the present paper, the literature has observed the significance of \cite[Theorem 1]{JinKobayashi2019} in the context of numerical approximation of time-changed SDEs, and many results that previously used \cite[Theorem 1]{JinKobayashi2019} for a limited class of inverse subordinators can be immediately generalized.

Let $D=(D_s)_{s\ge0}$ be a subordinator, that is, a one-dimensional nondecreasing
L\'evy process such that $D_0=0$, with Laplace exponent $\phi$. Thus, $\E[e^{-\lambda D_s}]=e^{-s\phi(\lambda)}$ for all $\lambda,s\ge0$.
Let $E$ be the right-continuous inverse of $D$, defined by
\[
E_t\coloneqq \inf\{s>0:D_s>t\},\qquad t\ge0.
\]
The process $E=(E_t)_{t\ge 0}$ is called an \textit{inverse subordinator} and will be used as a random time change throughout this paper. 
It is well known that Laplace exponents for subordinators are in one-to-one correspondence with Bernstein functions starting from zero, and that a given Laplace exponent $\phi$ has the representation
\begin{align}\label{bernstein-exp}
    \phi(\lambda)=\mathfrak{b}\lambda+\int_0^\infty\left(1-e^{-\lambda r}\right)\nu(\dd r),\qquad\lambda\ge0.
\end{align}
Here $\mathfrak{b}\ge0$ is called the drift coefficient, and $\nu$ is called the L\'evy measure and satisfies $\int_0^\infty (1\wedge r)\,\nu(\dd r)<\infty$. For further details, see \cite{Bertoin_subordinator,bernstein}. 
Notice that the property $\phi(\infty)=\infty$ is equivalent to the condition $\mathfrak{b}>0$ or $\nu(0,\infty)=\infty$, and from a probabilistic viewpoint, the latter is equivalent to the condition that sample paths of $D$ are strictly increasing. In other words, if $\phi(\infty)=\infty$, then $D$ is not a compound Poisson process, and consequently, sample paths of the inverse subordinator $E$ are continuous and nondecreasing. On the other hand, if $\phi(\infty)<\infty$, then $D$ is a compound Poisson process, which implies that sample paths of $E$ can increase by jumps only. 
In the discussion below, we will exclude the case when $\mathfrak{b}=0$ and $\nu=0$ simultaneously, which corresponds to the identically zero process with $\phi\equiv 0$.

A useful observation about a Bernstein function is that
\begin{align}\label{phi' bound}
    \phi'(\lambda)\le \frac{\phi(\lambda)}{\lambda},\quad\lambda>0, 
\end{align}
which follows from its integral representation. Moreover, it immediately follows that $\lambda\mapsto \phi(\lambda)/\lambda$ is decreasing.

In some of our results, we will need some finer properties of $\phi$. We say a function $f:[0,\infty)\to[0,\infty)$ satisfies the \textit{weak lower scaling condition (WLSC)} at infinity if there exist {$\alpha\in\R$,} $c_0\in(0,1]$, and $\lambda_0\ge 0$ such that 
\begin{align*}
	f(x \lambda)\ge c_0 x^\alpha f(\lambda), \quad x\ge 1,\, \lambda\ge\lambda_0.
\end{align*}
Following the usual notation, we write $f\in \textrm{WLSC}(\alpha,c_0,\lambda_0)$ in this case. 
Similarly, we say that $f$ satisfies the \textit{weak upper scaling condition (WUSC)} at infinity if there exist {$\alpha\in\R$,} $C_0\ge 1$, and $\lambda_0\ge 0$ such that 
 \begin{align*}
 	f(x \lambda)\le C_0 x^\alpha f (\lambda), \quad x\ge 1, \lambda\ge\lambda_0.
 \end{align*}
 We write $f\in \textrm{WUSC}(\alpha,C_0,\lambda_0)$ for short. 
Note that if $f$ is continuous, positive on $(0,\infty)$, and satisfies WLSC or WUSC with some $\lambda_0>0$, then the same type of inequality also holds with even smaller $\lambda_0>0$, with a possibly different constant $c_0$ or $C_0$. However, it may be impossible to extend the property to $\lambda_0=0$.

If $D$ is an $\alpha$-stable subordinator with $\alpha\in(0,1)$, then its Laplace exponent satisfies $\phi(\lambda)=\lambda^\alpha$, and consequently, $\phi\in \textrm{WLSC}(\alpha,1,0)\cap \textrm{WUSC}(\alpha,1,0)$; i.e., the scaling is exact and global.
Moreover, if $\phi$ is regularly varying at infinity with index $\alpha\in(0,1)$, then 
for any fixed $\varepsilon\in(0,\min(\alpha,1-\alpha))$, there exist constants $c_0\in(0,1]$, $C_0\ge 1$, and $\lambda_0\ge 0$ such that $\phi\in\textrm{WLSC}(\alpha-\varepsilon,c_0,\lambda_0)\cap \textrm{WUSC}(\alpha+\varepsilon,C_0,\lambda_0)$. 
Also note that if $\phi\in \textrm{WUSC}(\alpha,c_0,\lambda_0)$, then $\phi$ has no drift; i.e., $\mathfrak{b}=0$ in \eqref{bernstein-exp}. On the other hand, if $\phi\in\textrm{WLSC}(\alpha,c_0,\lambda_0)$, then $\phi(\infty)=\infty$, so the corresponding subordinator $D$ is not a compound Poisson process.
A number of examples of Laplace exponents $\phi$ satisfying \textrm{WLSC} or \textrm{WUSC} can be found in \cite[p. 409]{KSV_bhpinf-PA14} and \cite[p. 41]{KSV_heat-PA18}.

Some of our results assume that a given Laplace exponent $\phi$ is a \textit{special Bernstein function}, which means that $\lambda\mapsto \lambda/\phi(\lambda)$ is also a Bernstein function. Special Bernstein functions include all complete Bernstein functions, for which the L\'evy measure $\nu$  admits a completely monotone density. For details, see \cite{bernstein}, and in particular, \cite[Chapter 20]{bernstein} for a vast list of explicit complete Bernstein functions, which include Laplace exponents for stable subordinators,  exponentially tempered stable subordinators, and gamma subordinators.

\subsection{Incremental moments for inverse subordinators}

For a subordinator $D$, define the renewal function $U:[0,\infty)\to [0,\infty)$ by 
\[
U(t)\coloneqq \E\left[\int_0^\infty\1_{[0,t]}(D_s)ds\right],
\]
which measures how much time $D$ spends in the interval $[0,t]$. It is easy to see that
\begin{align}\label{U-def}
	U(t)=\E[E_t], \quad t\ge 0.
\end{align}
By \cite[Proposition 1 in Chapter III]{Bertoin_Levy}, $U(t)$ is comparable to $1/\phi(t^{-1})$ in the sense that there exist constants $c_1,C_1>0$ such that 
\begin{align}\label{U-and-phi}
	 \frac{c_1}{\phi(t^{-1})}\le U(t) \le \frac{C_1}{\phi(t^{-1})}, \quad t>0. 
\end{align}

The following is the first main theorem of this paper. It generalizes the known inequalities \eqref{Eincrements-stable} valid for inverse \textit{stable} subordinators to general inverse subordinators. Note that we follow the usual convention that $1/0\coloneqq \infty$ and $1/\infty\coloneqq 0$, and that $f^n$ denotes the $n$th power of a given function $f$. 

\begin{theorem}\label{Theorem_increments}
Let $E$ be the inverse of a subordinator $D$ with Laplace exponent $\phi$.
\begin{enumerate} 
\item [(a)]
{
For all $t\ge0$ and $p>0$, it holds that 
\begin{align}\label{E-moment}
    \max_{\varepsilon\in[0,1]}\eps^p\left(\frac{e^{-\eps}-e^{-1}}{1-e^{-1}}\right)\frac{1}{\phi^p(t^{-1})}\le \E[ E_t^p]\le \frac{e\hspace{1pt}\Gamma(p+1)}{\phi^p(t^{-1})}.
\end{align}
}
\item[(b)] (i)
For all $0\le a<b< \infty$ and $n\in \N$, it holds that 
\begin{align}\label{E-increment-general}
    \E[(E_b-E_a)^n]\le n! \left(U(b)-U(a)\right)U^{n-1}(b-a)\le \frac{e^{n-1}n!}{\phi^{n-1}((b-a)^{-1})}\left(U(b)-U(a)\right).
\end{align}
In particular, if $\phi$ is a special Bernstein function, i.e., $\lambda\mapsto \lambda/\phi(\lambda)$ is also a Bernstein function, then for all $0\le a<b<\infty$ and all $n\in\N$, it holds that
\begin{align}\label{E-increment-0}
    \E[(E_b-E_a)^n]\le \frac{e^{n-1}n!{(1-2e^{-1})^{-1}}}{\phi^{n-1}((b-a)^{-1})}\left(\frac{1}{\phi(b^{-1})}-\frac{1}{\phi(a^{-1})}\right).
\end{align} 
(ii) If $\phi$ is a special Bernstein function, i.e., $\lambda\mapsto \lambda/\phi(\lambda)$ is also a Bernstein function, and $\phi\in \textrm{WLSC}(\alpha_0,c_0,\lambda_0)\cap \textrm{WUSC}(\alpha_1,c_1,\lambda_0)$ with $0<\alpha_0\le \alpha_1<1$, then there exists a constant $C=C(n,\alpha_0,\alpha_1,c_0,c_1,\lambda_0)>0$ such that for all $0\le a<b<\lambda_0^{-1}$ and all $n\in\N$, it holds that
\begin{align}\label{E-increment-1}
    \E[(E_b-E_a)^n]\ge \frac{{Cn!}}{\phi^{n-1}((b-a)^{-1})}\left(\frac{1}{\phi(b^{-1})}-\frac{1}{\phi(a^{-1})}\right).
\end{align}
\end{enumerate}
\end{theorem}

\begin{proof}
(a) For any $\lambda>0$,
	\begin{align*}
		\E[ E_t^p]&=\int_0^\infty px^{p-1}\P(E_t>x)dx\le \int_0^\infty px^{p-1}\P(D_x\le t)dx= \int_0^\infty px^{p-1}\P(e^{-\lambda D_x}\ge e^{-\lambda t})dx\\
        &\le \int_0^\infty px^{p-1}e^{\lambda t}\E e^{-\lambda D_x}dx=pe^{\lambda t}\int_0^\infty x^{p-1}e^{-x \phi(\lambda)}dx=\frac{\Gamma(p+1)e^{\lambda t}}{\phi(\lambda)^p},
	\end{align*}
	where the inequality in the second line follows from Markov's inequality. The upper bound in \eqref{E-moment} now follows by taking $\lambda=t^{-1}$.  

{To derive the lower bound in \eqref{E-moment}, fix $\eps\in[0,1]$, and observe by Markov's inequality that 
\begin{align*}
	\phi^p(t^{-1})\E[E_t^p]
	\ge \eps^p\P\Bigl(E_t>\frac{\eps}{\phi(t^{-1})}\Bigr)
	\ge \eps^p\P\Bigl(D_{\frac{\eps}{\phi(t^{-1})}} < t\Bigr).
\end{align*}
For any $a>0$, 
\begin{align*}
	\P(D_a\ge t) 
	= \P(1-e^{-t^{-1}D_a} \ge 1-e^{-1})
	\le \frac{1-\E[e^{-t^{-1}D_a}]}{1-e^{-1}}
	= \frac{1-e^{-a\phi(t^{-1})}}{1-e^{-1}},
\end{align*}
and hence,
\[
	\P\Bigl(D_{\frac{\eps}{\phi(t^{-1})}} < t\Bigr)
	\ge 1- \frac{1-e^{-\eps}}{1-e^{-1}}
	=\frac{e^{-\eps}-e^{-1}}{1-e^{-1}}. 
\]
Putting these together gives
\[
	\phi^p(t^{-1})\E[E_t^p]
	\ge \max_{\eps\in[0,1]}\eps^p\left(\frac{e^{-\eps}-e^{-1}}{1-e^{-1}}\right),
\]
as desired.}

(b)    By applying \cite[Theorem 1]{Lageras}, we have
    \begin{align}\label{lageras}
        \E[(E_b-E_a)^n]=n!\int_{(a,b]}U(dx_1)\int_{[x_1,b]}U(dx_2-x_1)\cdots\int_{[x_{n-1},b]}U(dx_n-x_{n-1}).
    \end{align}
    Note that \cite[Theorem 1]{Lageras} was obtained by using Cox processes under the assumption that $\phi(\infty)=\infty$, i.e., compound Poisson processes were excluded. However, the statement was extended to include the case $\phi(\infty)<\infty$ via the use of the multi-point tail distribution function of the inverse subordinator in \cite[Theorem 4.2]{VeilletteTaqqu-letters}.
    Now, using the monotonicity of $U$ and combining \eqref{U-def} and \eqref{E-moment} gives
    \begin{align*}
        \int_{[x_{n-1},b]}U(dx_n-x_{n-1})=U(b-x_{n-1})\le U(b-a)\le \frac{e}{\phi((b-a)^{-1})},
    \end{align*}
    By iterating this argument, \eqref{lageras} gives
    \begin{align*}
        E[(E_b-E_a)^n]&\le n!U^{n-1}(b-a)\int_{(a,b]}U(dx_1)=n!U^{n-1}(b-a)\big(U(b)-U(a)\big)\\&\le n!\frac{e^{n-1}}{\phi^{n-1}((b-a)^{-1})}\left(U(b)-U(a)\right),
    \end{align*}
    thereby establishing \eqref{E-increment-general}.
    
    Let us now assume that $\phi$ is a special Bernstein function. Then by \cite[Chapter 13]{bernstein}, $U$ is absolutely continuous on $(0,\infty)$ with a decreasing density $u$, possibly with an atom at 0. More precisely, it holds that $U(t)-U(0)=\int_0^tu(s)\,ds$ and $U(0)=\frac{1}{\phi(\infty)}$. In other words, $U(ds)=\1_{(0,\infty)}(s)u(s)ds+\frac{1}{\phi(\infty)}\delta_0(ds)$; see \cite[Theorem 13.3 \& p. 198]{bernstein}.
    Therefore, by \cite[p. 74]{Bertoin_Levy},    
    \begin{align}
    \begin{split}\label{U-Laplace}
        \frac{1}{\phi(\infty)}+\int_0^\infty e^{-\lambda s}u(s)\,ds&=\int_{[0,\infty)} e^{-\lambda s}U(ds)=\E\left[\int_0^\infty e^{-\lambda D_s} ds\right]\\
        &=\int_0^\infty \E\left[e^{-\lambda D_s} \right]ds=\int_0^\infty e^{-s\phi(\lambda)}ds=\frac{1}{\phi(\lambda)}.       
    \end{split}
    \end{align}
    Hence, {application of \cite[Lemma A.1]{KimMimica2012} with $f\coloneqq 1/\phi$ and $\nu\coloneqq u$ yields}
    \begin{align}\label{u bound}
        u(t)\le {(1-2e^{-1})^{-1}}\frac{\phi'(t^{-1})}{t^2\phi^2(t^{-1})},\quad t>0.
    \end{align} 
By the change of variables $\phi(t^{-1})=h$, 
    \begin{align*}
        \int_a^b\frac{\phi'(t^{-1})}{t^2\phi^2(t^{-1})}dt
        =\int_{\phi(b^{-1})}^{\phi(a^{-1})}\frac{dh}{h^2}=\frac{1}{\phi(b^{-1})}-\frac{1}{\phi(a^{-1})},
    \end{align*}
    and hence,
    \begin{equation*}
        U(b)-U(a)=\int_a^b u(s)\,ds\le{(1-2e^{-1})^{-1}}\left( \frac{1}{\phi(b^{-1})}-\frac{1}{\phi(a^{-1})}\right).
    \end{equation*} 
    Inserting this into \eqref{E-increment-general} yields \eqref{E-increment-0}.
    
    We now turn our attention to the lower bound in \eqref{E-increment-1}, where it is assumed that $\phi$ is special, and $\phi\in \textrm{WLSC}(\alpha_0,c_0,\lambda_0)\cap \textrm{WUSC}(\alpha_1,c_1,\lambda_0)$ with $0<\alpha_0\le \alpha_1<1$, and $0\le a<b<\lambda_0^{-1}$. 
    Then $\mathfrak{b}=0$ and $\phi(\infty)=\infty$, and hence, $U(0)=0$ and $U(t)=\int_0^t u(s)ds$.
     First, by \cite[Lemma 1.3]{KimMimica2014} together with the assumption $\phi\in \textrm{WLSC}(\alpha_0,c_0,\lambda_0)$, there exists a constant $c_2>0$ such that
     \begin{align}\label{404}
        \phi'(\lambda)\ge c_2\frac{\phi(\lambda)}{\lambda}, \quad \lambda\ge \lambda_0.
     \end{align}
    Note that the reverse direction always holds for all $\lambda>0$, as stated in \eqref{phi' bound}.
    Next, it follows from \eqref{phi' bound}, \eqref{404}, and the assumption $\phi\in \textrm{WUSC}(\alpha_1,c_1,\lambda_0)$ that for all $x\ge 1$ and $\lambda\ge \lambda_0$,  
    \[
        \phi'(x\lambda)\le \frac{\phi(x\lambda)}{x\lambda}
        \le \frac{c_1 x^{\alpha_1}\phi(\lambda)}{x\lambda}
        \le \frac{c_1}{c_2}x^{\alpha_1-1}\phi'(\lambda).
    \]
    Therefore, $\phi'\in \textrm{WUSC}(\alpha_1-1,c_1/c_2,\lambda_0)$, and in particular, Assumption (A-3) in the paper \cite{KimMimica2012} is satisfied, which implies that the function $f\coloneqq 1/\phi$ satisfies Equation (A.3) of the same paper, as seen in the proof of \cite[Proposition 3.4]{KimMimica2012}. Then by \cite[Lemma A.2]{KimMimica2012}, not only the inequality in \eqref{u bound} but also the following reverse inequality holds with some constant $c_3>0$: 
    \begin{align}\label{u bound-lower}
        u(t)\ge c_3\frac{\phi'(t^{-1})}{t^2\phi^2(t^{-1})},\quad t\in(0,\lambda_0^{-1}].
    \end{align}
    Combining \eqref{404} and \eqref{u bound-lower} gives
    $u(t)\ge c_3c_2/(t\phi(t^{-1}))$ for all $t\in(0,\lambda_0^{-1}]$. 
    It then follows from the iterated integral representation \eqref{lageras}, using $U(dt)=u(t)dt$, where the variables satisfy $a\le x_1<x_2<\cdots<x_n\le b<\lambda_0^{-1}$, that
    \begin{align*}
        \E[(E_b-E_a)^n]&\ge n!(c_3c_2)^n\int_a^b\frac{1}{x_1\phi(x_1^{-1})}\biggl(\int_{x_1}^b\cdots \biggl(\int_{x_{n-1}}^b \frac{1}{(x_n-x_{n-1})\phi((x_n-x_{n-1})^{-1})}\, dx_n\biggr)\cdots dx_2\biggr)\,dx_1.
\end{align*}
By the change of variables $(x_1,x_2,\cdots,x_n)\leftrightarrow (z_1,z_2,\cdots,z_n)$ given by 
\begin{align*}
	x_1&=a+(b-a)z_1, \\
	x_2&=x_1+(b-a)(1-z_1)z_2,\\
	x_3&=x_2+(b-a)(1-z_1)(1-z_2)z_3,\\
	& \ \, \vdots\\
	x_n&=x_{n-1}+(b-a)(1-z_1)(1-z_2)\cdots (1-z_{n-1})z_n,
\end{align*}
we obtain
\begin{align*}
\E[(E_b-E_a)^n]&\ge  n!(c_3c_2)^n(b-a)\int_0^1 \frac{1}{f(z_1)} \biggl(\int_0^1 \frac{1}{g_2(z_1,z_2)} \biggl(\int_0^1 \frac{1}{g_3(z_1,z_2,z_3)}\\
	&\ \ \ \ \ \cdots \biggl(\int_0^1 \frac{1}{g_n(z_1,\cdots,z_n)}\,dz_n\biggr)\,\cdots\,dz_3\biggr)\,dz_2\biggr)\,dz_1,
\end{align*}
where the functions $f,g_2,\ldots,g_n$ are defined by
\begin{align*}
	f(z_1)\coloneqq [a+(b-a)z_1]\phi\biggl(\frac{1}{a+(b-a)z_1}\biggr), 
\end{align*}
and
\begin{align*}
	g_k(z_1,\ldots,z_k)&\coloneqq 
	z_k\phi\biggl(\frac{1}{(b-a)(1-z_1)(1-z_2)\cdots (1-z_{k-1})z_k}\biggr).
\end{align*}
Using $bz_1\le a+(b-a)z_1 \le b$ and $\phi\in \textrm{WUSC}(\alpha_1,c_1,\lambda_0)$ with $x=1/z_1\ge 1$ and $\lambda=1/b$ gives 
\[
	f(z_1)=[a+(b-a)z_1]\phi\biggl(\frac{1}{a+(b-a)z_1}\biggr)
	 \le b\phi\biggl(\frac{1}{bz_1}\biggr)
	 \le c_1 b  z_1^{-\alpha_1}\phi\biggl(\frac1b\biggr). 
\]
On the other hand, for each $2\le k\le n$, by using $\phi\in \textrm{WUSC}(\alpha_1,c_1,\lambda_0)$ with $x=1/[(1-z_1)(1-z_2)\cdots (1-z_{k-1})z_k]\ge 1$ and $\lambda=1/(b-a)> 1/b$, we obtain 
\begin{align*}
	g_k(z_1,\ldots,z_k)
	 &\le c_1z_k^{1-\alpha_1} [(1-z_1)(1-z_2)\cdots (1-z_{k-1})]^{-\alpha_1}  \phi\biggl(\frac{1}{b-a}\biggr). 
\end{align*}
Therefore, setting $c_4=c_3c_2/c_1$ gives
\begin{align*}
	&\E\left[ (E_b-E_a)^n\right]\\
	&\ge \frac{n!c_4^n}{b \phi(b^{-1})}\frac{b-a}{\phi^{n-1}((b-a)^{-1})}\int_0^1 z_1^{\alpha_1} \biggl(\int_0^1(1-z_1)^{\alpha_1} z_2^{\alpha_1-1} \biggl(\int_0^1 (1-z_1)^{\alpha_1}(1-z_2)^{\alpha_1} z_3^{\alpha_1-1}\\
	&\ \ \ \ \ \cdots \biggl(\int_0^1 (1-z_1)^{\alpha_1}(1-z_2)^{\alpha_1}\cdots (1-z_{n-1})^{\alpha_1} z_n^{\alpha_1-1} \, dz_n\biggr)\,\cdots \,dz_3\biggr)\,dz_2\biggr)\,dz_1. 
\end{align*}
{For $n\ge 2$,} the latter multiple integral simplifies to
\begin{align*} 
	I(\alpha_1,n) \coloneqq &\int_0^1 z_n^{\alpha_1-1} \, dz_n\times \int_0^1 z_1^{\alpha_1}(1-z_1)^{(n-1)\alpha_1}\,dz_1 \times\prod_{\ell=2}^{n-1}\int_0^1 z_\ell^{\alpha_1-1}(1-z_\ell)^{(n-\ell)\alpha_1}\,dz_\ell \\
	= &\frac 1{\alpha_1} B(1+\alpha_1,1+(n-1)\alpha_1)\prod_{k=1}^{n-2} B(\alpha_1,1+k\alpha_1)
    =\frac{\Gamma^2(1+\alpha_1)\Gamma^{n-2}(\alpha_1)}{\alpha_1 \Gamma(2+n\alpha_1)}
    =\frac{\alpha_1\Gamma^{n}(\alpha_1)}{ \Gamma(2+n\alpha_1)},
\end{align*}
{while $I(\alpha_1,1)=\frac{1}{1+\alpha_1}$.}
Therefore, 
\begin{align}\label{406}
    &\E[(E_b-E_a)^n]\ge {n!c_4^nI(\alpha_1,n)}\frac{b-a}{b\phi(b^{-1})}\frac{1}{\phi^{n-1}((b-a)^{-1})}.
\end{align}
Finally, $\lambda\mapsto \phi(\lambda)/\lambda$ is decreasing due to \eqref{phi' bound}, and hence, $a\phi(a^{-1})\le b\phi(b^{-1})$. This is equivalent to $1-\phi(b^{-1})/\phi(a^{-1})
\le 1-a/b$, which is further equivalent to 
\begin{align}\label{407}
\frac{1}{\phi(b^{-1})}-\frac{1}{\phi(a^{-1})}
 \le \frac{b-a}{b\,\phi(b^{-1})}.
\end{align}
Combining \eqref{406} and \eqref{407} yields the lower bound in \eqref{E-increment-1}. 
\end{proof}
\begin{remark}\label{r:incremental_moments}
\begin{enumerate}[(a)]
    \item Theorem \ref{Theorem_increments}(b) shows that for a wide class of Laplace exponents $\phi$, the incremental moments of the inverse subordinators behave as follows:
    \[
        \E[(E_b-E_a)^n]\asymp\frac{1}{\phi^{n-1}((b-a)^{-1})}\left(\frac{1}{\phi(b^{-1})}-\frac{1}{\phi(a^{-1})}\right).
    \]
    \item {Regarding Theorem \ref{Theorem_increments}(a), if $p\ge 1$, one could easily obtain a lower bound of $\E[ E_t^p]$ by combining \eqref{U-and-phi} with Jensen's inequality as follows:
	\begin{align*}
		\frac{c_1^p}{\phi^p(t^{-1})}\le (\E [E_t])^p\le  \E[ E_t^p].
	\end{align*}
    The proof of \cite[Proposition 1 in Chapter III]{Bertoin_Levy} indicates that we may take $c_1=\frac{1}{4\log 8}$. However, this argument is not valid for $0<p<1$, and even for $p\ge 1$, this particular lower bound $\frac{c_1^p}{\phi^p(t^{-1})}$ is not as good as the one given in \eqref{E-moment}. Indeed, considering $\eps=\frac 12$ gives
\[
    \max_{\eps\in[0,1]}\eps^p\left(\frac{e^{-\eps}-e^{-1}}{1-e^{-1}}\right)
    \ge 2^{-p}\left(\frac{e^{-1/2}-e^{-1}}{1-e^{-1}}\right)
    =\frac{2^{-p}}{\sqrt{e}+1}
    >\frac{1}{(4\log 8)^p}.
\]
}

\end{enumerate}
\end{remark}

\subsection{Exponential moments of powers of inverse subordinators}

Given an inverse subordinator $E$ with $\phi(\infty)=\infty$,  it holds that $\mathbb{E}[e^{cE_t}]<\infty$ for any constants $c,t>0$ (see, e.g., \cite[Lemma 3.4]{JumKobayashi} and \cite[Lemma A.3]{BCT2026}). 
In particular, it follows from \cite[Lemma A.3]{BCT2026} that 
\begin{align}\label{A.3}
    \mathbb{E}[e^{cE_t}]\le \tilde{c}\hspace{1pt}e^{\tilde{c}t}
\end{align}
for some constant $\tilde{c}>0$ depending on $c$ and $\phi$.
It has also been proved in \cite[Theorem 1]{JinKobayashi2019} that if the Laplace exponent $\phi$ of the underlying subordinator $D$ is regularly varying with index $\alpha\in(0,1)$, then 
\[
	\mathbb{E}[e^{cE_t^r}]<\infty \ \ \textrm{if} \ 0<r<\frac{1}{1-\alpha},
\]
and
\[
	\mathbb{E}[e^{cE_t^r}]=\infty \ \ \textrm{if} \ r>\frac{1}{1-\alpha}. 
\]
Hence, for the expectation $\mathbb{E}[e^{cE_t^2}]$ (with $r=2$), it follows that $\mathbb{E}[e^{cE_t^2}]<\infty$ if $\alpha>1/2$, while $\mathbb{E}[e^{cE_t^2}]=\infty$ if $\alpha<1/2$. 
This fact is easily verified when the subordinator $D$ is $\alpha$-stable; see \cite[Remark 6(2)]{JinKobayashi2019}.

The criterion for the finiteness of $\mathbb{E}[e^{cE_t^r}]$ has played a crucial role in the error analysis of the numerical approximation of time-changed SDEs; see Remarks \ref{rem:expo_moment} and \ref{rem:Zuo}. 
The following statement generalizes \cite[Theorem 1]{JinKobayashi2019} to situations where the Laplace exponent $\phi$ in \eqref{bernstein-exp} satisfies a weak scaling condition only.

\begin{theorem}\label{theorem_expo_moment}
Let $E$ be the inverse of a subordinator $D$ with Laplace exponent $\phi$. 
Let $c,t>0$ be constants. 
\begin{enumerate}[(a)]
\item If $\nu(0,\infty)=\infty$ and $\phi$ satisfies
\begin{align}\label{phi-doubleprime}
	-\phi'' \in \textrm{WLSC}(\alpha-2,\tilde{c}_0,\lambda_0),
\end{align}
 where $\alpha\in(0,1)$, $\tilde{c}_0\in(0,1]$, and $\lambda_0\ge 0$, then $\mathbb{E}[e^{cE_t^r}]<\infty$ for  any $r\in(0,\frac{1}{1-\alpha})$.
\item If $\phi$ satisfies
\begin{align}\label{phi-USC}
	\phi \in \textrm{WUSC}(\alpha,C_0,\lambda_0)
\end{align}
 where $\alpha\in(0,1)$, $C_0\ge 1$, and $\lambda_0\ge 0$, then $\mathbb{E}[e^{cE_t^r}]=\infty$ for any $r>\frac{1}{1-\alpha}$. 
 \end{enumerate}
\end{theorem}

\begin{proof}
(a) By \cite[Propositions 2.2 and 2.3]{GLT2023} (also see Remark 2.5 of that paper), condition \eqref{phi-doubleprime} implies that 
\begin{align}\label{phi-prime}
	\phi' \in \textrm{WLSC}(\alpha-1,\tilde{c}_0,\lambda_0).
\end{align}
As in \cite{JinKobayashi2019}, consider the following auxiliary functions initially introduced in \cite{JainPruitt}: 
\[
	g(s)=\phi'(s) \ \  \textrm{and} \ \ R(s)=\phi(s)-s\phi'(s).
\]
Let
\[
	x(s)=g(s^r)=\phi'(s^r). 
\]
Then for any $w\ge 1$ and $\kappa> \lambda_0^{1/r}$, since $w^r\ge 1$ and $\kappa^r> \lambda_0$, the lower scaling condition \eqref{phi-prime} gives
\[
	x(w \kappa)=\phi'(w^r \kappa^r)
	\ge \tilde{c}_0 (w^r)^{\alpha-1} \phi'(\kappa^r)
	= \tilde{c}_0 w^{r(\alpha-1)} x(\kappa),  
\]
and hence, $x\in \textrm{WLSC}(r(\alpha-1),\tilde{c}_0,\lambda_0^{1/r})$. 
This immediately implies that the function $z(s)=sx(s)$ satisfies 
\[
	z\in \textrm{WLSC}(r(\alpha-1)+1,\tilde{c}_0,\lambda_0^{1/r}).
\] 
In particular, for a fixed $\kappa>\lambda_0^{1/r}$, since $r(\alpha-1)+1>0$ if $r\in(0,\frac{1}{1-\alpha})$, 
\[
	\liminf_{s\to \infty}z(s) =\liminf_{w\to \infty}z(w\kappa) 
	\ge \liminf_{w\to \infty}\tilde{c}_0w^{r(\alpha-1)+1} z(\kappa)=\infty.
\]
Here, we used $z(\kappa)>0$. To verify the latter, suppose $z(\kappa)=0$. Then $\phi'(\kappa^r)=x(\kappa)=0$. But since $\phi''\le 0$ and $\phi'\ge 0$, it follows that $\phi'=0$ on $[\kappa^r,\infty)$, and hence, $\phi(s)$ is constant on $[\kappa^r,\infty)$, which contradicts the fact that $\phi$ is strictly increasing on $(0,\infty)$.  (In fact, for $s_1<s_2$, we can see that 
$
	\phi(s_2)-\phi(s_1)
	 \ge \int_0^\infty (1-e^{-s_2 y})\,\nu(dy) - \int_0^\infty (1-e^{-s_1 y})\,\nu(dy)
	 = \int_0^\infty (e^{-s_1 y}-e^{-s_2 y})\,\nu(dy)>0
$.)

On the other hand, observe that 
\[
	R'(s)=\phi'(s)-(s\phi'(s))' = s(-\phi''(s)).
\]
For any $w\ge 1$ and $\lambda> \lambda_0$, by \eqref{phi-doubleprime}, 
\[
	R'(w\lambda) = w\lambda(-\phi''(w\lambda))
	\ge w\lambda \cdot \tilde{c}_0w^{\alpha-2}(-\phi''(\lambda))
	=\tilde{c}_0w^{\alpha-1}R'(\lambda).
\]
Therefore, $R'\in \textrm{WLSC}(\alpha-1,\tilde{c}_0,\lambda_0)$. 
Now, for fixed $w\ge 1$ and $y> \lambda_0$, 
\begin{align*}
	R(w y)-R(w \lambda_0) 
	&=\int_{w \lambda_0}^{w y} R'(s)\,\dd s 
	=w\int_{\lambda_0}^{y} R'(w u)\,du \\
	&\ge w\int_{\lambda_0}^{y} \tilde{c}_0w^{\alpha-1}R'(u)\,du
	=\tilde{c}_0w^\alpha \{R(y)-R(\lambda_0)\}, 
\end{align*}
from which it follows that 
\begin{align}\label{R-function}
	R(w y)\ge R(w \lambda_0) + \tilde{c}_0w^\alpha \{R(y)-R(\lambda_0)\}.
\end{align}
Recall that $R$ is a continuous, nondecreasing function such that $R(0)=0$ and $R(\infty)=\nu(0,\infty)$. Since $\nu(0,\infty)=\infty$ by assumption,
we can take $y> \lambda_0$ large enough so that $R(y)>R(\lambda_0)$. Then by \eqref{R-function} and the inequality $\alpha-\frac{r-1}{r}>0$, which is true since $r\in(0,\frac{1}{1-\alpha})$, 
\begin{align*}
	\liminf_{s\to \infty} \frac{R(s^r)}{s^{r-1}}
	&=\liminf_{w\to \infty} \frac{R(((w y)^{1/r})^r)}{((w y)^{1/r})^{r-1}}
	= \liminf_{w\to \infty} \frac{R(w y)}{(w y)^{(r-1)/r}}\\
	&\ge  \liminf_{w\to \infty} \frac{R(w \lambda_0) + \tilde{c}_0w^\alpha \{R(y)-R(\lambda_0)\}}{(w y)^{(r-1)/r}}\\
	&\ge  \liminf_{w\to \infty} \frac{\tilde{c}_0 \{R(y)-R(\lambda_0)\}}{y^{(r-1)/r}} w^{\alpha-(r-1)/r}=\infty. 
\end{align*}
Thus, we have verified that $z(s)=sx(s)\to \infty$ and $\frac{R(s^r)}{s^{r-1}}\to \infty$ as $s\to\infty$, and application of \cite[Corollary 3 (1)]{JinKobayashi2019} yields the finiteness of $\mathbb{E}[e^{cE_t^r}]$.
(Notice that, in \cite{JinKobayashi2019}, it is assumed that $\mathfrak{b}=0$ and $\nu(0,\infty)=\infty$; however, the proofs of \cite[Proposition 2, Corollary 3]{JinKobayashi2019} show that the statements are actually valid for all $\mathfrak{b}$ and $\nu$.) 

(b) By \cite[Corollary 3 (2)]{JinKobayashi2019}, we only need to prove that $z(s)\to 0$ and $\frac{R(s^r)}{s^{r-1}}\to 0$ as $s\to\infty$. To this end, it suffices to show that  
\[
	\lim_{s\to\infty} \frac{\phi(s^r)}{s^{r-1}}=0
\]
since $0\le s^r\phi'(s^r)\le \phi(s^r)$ as $R(s^r)\ge 0$. 
Now, by the weak upper scaling condition \eqref{phi-USC}, for any fixed $\lambda_1>\lambda_0\ge 0$, 
\begin{align*}
	\limsup_{s\to\infty}\frac{\phi(s^r)}{s^{r-1}} 
	&=\limsup_{w\to \infty}\frac{\phi(w \lambda_1)}{(w \lambda_1)^{(r-1)/r}}
	\le \limsup_{w\to \infty}\frac{C_0w^\alpha\phi(\lambda_1)}{(w \lambda_1)^{(r-1)/r}}\\
	&= \limsup_{w\to \infty}\frac{C_0w^{\alpha-(r-1)/r}\phi(\lambda_1)}{\lambda_1^{(r-1)/r}}=0, 
\end{align*}
where the last equality follows from the assumption $r>\frac{1}{1-\alpha}$, which is equivalent to $\alpha-\frac{r-1}{r}<0$.
\end{proof}

\begin{remark}
\begin{enumerate}[(a)]
    \item 
    By \cite[Proposition 2.3]{GLT2023}, condition \eqref{phi-prime} for some $\tilde{c}_0$ and the following condition for some $c_0$ are equivalent:
\begin{align}\label{phi-0}
	\phi \in \textrm{WLSC}(\alpha,c_0,\lambda_0).
\end{align}
 However, even though \eqref{phi-doubleprime} implies \eqref{phi-prime}, the converse may not hold.  The above proof for part (a) relies on \eqref{phi-doubleprime}, and hence, it is not sufficient to assume \eqref{phi-prime} or \eqref{phi-0}.

\item Conditions \eqref{phi-doubleprime} and \eqref{phi-prime} become equivalent if $\phi$ satisfies both weak upper and lower scaling conditions; see \cite[Propositions 2.2 and 2.6, Corollary 2.7]{GLT2023}. 
\end{enumerate}
\end{remark}

\begin{remark}
\label{rem:expo_moment}
Although the criteria for the existence and non-existence of $\E[e^{cE_t^r}]$ established in Theorem \ref{theorem_expo_moment} are not used in the present paper, they may be of independent interest to some readers. In fact, they are generally useful when analyzing the sizes of the errors associated with numerical approximation schemes of time-changed SDEs, including those of the form \eqref{403} and \eqref{402}. 
In fact, \cite[Theorem 1]{JinKobayashi2019}, which is a less general version of Theorem \ref{theorem_expo_moment}, has been used in the literature; see e.g., \cite{JinKobayashi2019,JinKobayashi2021,Zuo2026}.
A further discussion on Theorem \ref{theorem_expo_moment} will be provided in Remark \ref{rem:Zuo}.
\end{remark}

\section{Strong approximation of time-changed SDEs with piecewise-constant interpolation: the \textit{asynchronous} case}\label{s:approx-1}

Let $(\Omega,\F,\P)$ be a probability space equipped with a filtration $(\F_t)_{t\ge 0}$ satisfying the usual conditions. 
Let $B$ be an $m$-dimensional $(\F_t)$-adapted Brownian motion which is independent of an $(\F_t)$-adapted subordinator $D$. 
Since the Brownian motion $B$ and the subordinator $D$ are assumed independent, it is possible to set up $B$ and $D$ on a product space $\Omega=\Omega_B\times \Omega_D$ with product measure $\P=\P_B\times \P_D$ with obvious notation. We use this setting in the proofs of all the remaining statements. Let $\E_B$, $\E_D$ and $\E$ denote the expectations under the probability measures $\P_B$, $\P_D$ and $\P$, respectively.

Throughout this section, we assume that the Laplace exponent $\phi$ of the subordinator $D$ satisfies $\mathfrak b>0$ or $\nu(0,\infty)=\infty$ in representation \eqref{bernstein-exp}. Recall that this is equivalent to assuming $\phi(\infty)=\infty$, or, in other words, $D$ is strictly increasing and is not a compound Poisson process. Consequently, the inverse $E$ of $D$ has continuous paths. 

Fix $t_0\ge 0$ and consider the time-changed SDE
\begin{align}\label{SDE_001-a}	
	X_t=x_{t_0}+\int_{t_0}^t F(s,X_s)\, \dd E_s+\int_{t_0}^t G(s,X_s)\,\dd B_{E_s} \ \ \textrm{for} \ t\in[t_0,T],
\end{align}
where $x_{t_0}\in\mathbb{R}^d$ is a non-random constant, $T>0$ is a fixed time horizon, and $F:[t_0,T]\times \mathbb{R}^d\to \mathbb{R}^d$ and $G:[t_0,T]\times \mathbb{R}^d\to \mathbb{R}^{d\times m}$ are measurable functions for which there exist constants $K>0$ and $\rho \in(0,1]$ such that
\begin{align}
	\label{SDE_condition1-a} &|F(t,x)-F(t,y)|+|G(t,x)-G(t,y)|\le K|x-y|,\\
	\label{SDE_condition2-a} &|F(t,x)|+|G(t,x)|\le K(1+|x|),\\
	\label{SDE_condition4-a} &|F(s,x)-F(t,x)|+|G(s,x)-G(t,x)|\le K(1+|x|)|s-t|^\rho
\end{align}  
for all $x,y\in\mathbb{R}^d$ and $s,t\in [t_0,T]$, with $|\cdot|$ denoting the Euclidean norms of appropriate dimensions. 
SDE \eqref{SDE_001-a} is understood within the framework of stochastic integrals driven by semimartingales, and 
the assumptions \eqref{SDE_condition1-a}--\eqref{SDE_condition2-a} guarantee the existence of a unique strong solution $X$ to SDE \eqref{SDE_001-a}; see the discussions in \cite{JinKobayashi2019,Kobayashi2011}. 
We say that the coefficients $F(s,X_s)$ and $G(s,X_s)$ of SDE \eqref{SDE_001-a} are \textit{asynchronous} because they do not share the same constant periods with the driving processes $E_s$ and $B_{E_s}$.
\textit{In the remainder of the paper, we assume that $m=d=1$ for simplicity of discussions and expressions;} an extension to a multidimensional case is straightforward.  

Following \cite{BCT2026}, we sample values of an inverse subordinator $E$ at \textit{deterministic} time points 
\begin{align}\label{deterministic-time-a}
 0\le t_0<t_1<t_2<\cdots<t_\ell<\cdots<t_L=T \ \  \textrm{with} \ \  {\Delta\coloneqq \max_{\ell\in\{0,1,2,\ldots,L-1\}}(t_{\ell+1}-t_\ell).}
 \end{align}
For each $t\in[t_0,T]$, let
\begin{align}\label{1505}
	\ell_t\coloneqq\max\{\ell\in \{0,1,2,\ldots,L\}: t_\ell\le t\},
\end{align}
so that 
\begin{align*}
	t_{\ell_t}\le t<t_{\ell_t+1} \ \ \textrm{for all} \ t\in[t_0,T), \ \ \textrm{and} \ \ t_{\ell_T}= T.
\end{align*}
Define a discrete-time process $(X^\Delta_{t_\ell})_{\ell\in \{0,1,2,\ldots,L\}}$ by setting 
\begin{align}
	&X^\Delta_{t_0}
	=x_{t_0}\label{Approx_001-a},\\
	&X^\Delta_{t_{\ell+1}}\!\!
	=X^\Delta_{t_{\ell}}+F(t_{\ell},X^\Delta_{t_{\ell}})(E_{t_{\ell+1}}-E_{t_{\ell}})+G(t_{\ell},X^\Delta_{t_{\ell}})(B_{E_{t_{\ell+1}}}-B_{E_{t_{\ell}}})
	\label{Approx_002-a}
\end{align}
for $\ell=0,1,2,\ldots,L-1$. 
By the assumption of independence between $B$ and $D$ (and hence between $B$ and $E$ as well), one can approximate the inverse subordinator $E$ first over the time points in \eqref{deterministic-time-a} and subsequently approximate the Brownian motion $B$ over the time steps $\{E_{t_0}, E_{t_1}, E_{t_2}, E_{t_3}, \ldots, E_{t_L}\}$.

In numerical analysis of SDEs, it is common to consider \textit{continuous interpolation}, where the path of an approximated solution is interpolated via the SDE. While this works well for theoretical analysis of errors ascribed to the approximation, it is impractical since one cannot actually visualize the entire sample path of the approximated solution. 
A remedy for this, with practical situations in mind, is to use an interpolation method that does not rely on the SDE itself. Here we consider the simplest of those: \textit{piecewise-constant interpolation}. Namely, given the approximate solution in discrete-time $(X^\Delta_{t_\ell})_{\ell\in \{0,1,2,\ldots,L\}}$ in \eqref{Approx_001-a}--\eqref{Approx_002-a}, define a continuous-time process $X^\Delta=(X^\Delta_t)_{t\in[t_0,T]}$ by 
\begin{align}
	X^\Delta_t\coloneqq X^\Delta_{t_{\ell_t}} \ \ \textrm{for any} \ t\in[t_0,T].\label{Approx_003-a}
\end{align}
The path of this approximated solution is visualizable on the entire time interval $[t_0,T]$. For the rest of the section, we will analyze the error associated with this particular approximation.  

Recall that 
an approximation process $X^\Delta$ with step size $\Delta>0$ is said to \textit{converge strongly to the solution $X$ uniformly over $[t_0,T]$ with order $\eta>0$} 
if there exist finite positive constants $C$ and $\Delta_0$ such that 
\[
	\mathbb{E}\left[\sup_{t_0\le t\le T}|X_t-X^\Delta_t|\right]\le C\Delta^\eta
\]
for all $\Delta\in(0,\Delta_0)$.

The following theorem establishes a uniform strong error estimate associated with the approximated solution $X^\Delta$ constructed above. 

\begin{theorem}\label{thm:strong-approx}
Let $E$ be the inverse of a subordinator $D$ with Laplace exponent $\phi$ in \eqref{bernstein-exp} with $\mathfrak{b}>0$ or $\nu(0,\infty)=\infty$. 
Let $B$ be a Brownian motion independent of $D$.
Let $X$ be the unique solution of SDE \eqref{SDE_001-a}, where the coefficients $F(t,x)$ and $G(t,x)$ satisfy conditions \eqref{SDE_condition1-a}--\eqref{SDE_condition4-a}. 
Let $X^\Delta$ be the approximation process constructed in \eqref{Approx_001-a}--\eqref{Approx_003-a} based on deterministic time points \eqref{deterministic-time-a} with $L\ge 2$.
Then, for every $\varepsilon\in(0,1)$, there exists a constant $C_1=C_1(x_{t_0},K,\phi,\varepsilon)>0$, independent of $\Delta$ and $T$, such that
\begin{align}\label{1427}
\E\left[\sup_{t_0\le t\le T}|X_t-X_t^\Delta|^2\right]\le C_1e^{C_1T}\left(\Delta^{2\rho}+\frac{1}{\phi(1/\Delta)^{1-\varepsilon}}\right),\qquad \Delta\in(0,1).
\end{align}
Moreover, if there exist constants $\alpha\in(0,1]$ and $c,\lambda_0>0$ such that $\phi(\lambda)\ge c\lambda^\alpha$ for all $\lambda\ge\lambda_0$, then 
there exists a constant $C_2=C_2(x_{t_0},K,\phi,\varepsilon)>0$, independent of $\Delta$ and $T$, such that
\begin{align}\label{1509}
\E\left[\sup_{t_0\le t\le T}|X_t-X_t^\Delta|^2\right]\le C_2e^{C_2T}\Delta^{\min\{2\rho,\alpha(1-\varepsilon)\}}, \qquad \Delta\in(0,1).
\end{align}
Consequently, $X^\Delta$ converges strongly to $X$ uniformly over $[t_0,T]$ with order $\eta=\min\left\{\rho,\frac{\alpha}{2}(1-\varepsilon)\right\}$.
\end{theorem}

The proof of Theorem \ref{thm:strong-approx} requires the following lemma. 

\begin{lemma}\label{lemma_Holder-a}
Let $X$ be as in Theorem \ref{thm:strong-approx}. Then for any $p\in[1,\infty)$, there exists a positive constant $C_{p,K}$ that depends on $p$ and $K$ such that 
\[
\E[Y^{(p)}_T]\leq2(1+|x_{t_0}|^p)\E\left[e^{C_{p,K}E_T}\right], 
\]
 where $Y^{(p)}_T\coloneqq 1+\sup_{t_0\le s\le T} |X_s|^p$. 
\end{lemma}
\begin{proof}
This statement appears in \cite[Lemma 8]{JinKobayashi2019} under the assumption that $\nu(0,\infty)=\infty$. However, that assumption is used only to guarantee that $E$ has continuous paths. Consequently, the same proof holds even when $\mathfrak b>0$ in \eqref{bernstein-exp}.
\end{proof}

\begin{proof}[{Proof of Theorem \ref{thm:strong-approx}}]

{In the proof, for any parameter $\theta$, the notation $C_\theta$ represents a positive constant depending on $\theta$, which may vary from line to line.}

Let $t\in[t_0,T]$. We first write
\begin{align}\label{Eq021-a}
	|X_t-X_t^\Delta|\le I_1(t)+I_2(t)+I_3(t)+I_4(t),
\end{align}
where
\begin{align*}
	&I_1(t)\coloneqq \left| \int_{t_0}^{t_{\ell_t}} (F(r,X_r)\hspace{-1pt}-F(t_{\ell_r},X^\Delta_r))\,\dd E_r \right|,
	&&I_2(t)\coloneqq \left| \int_{t_{\ell_t}}^t F(r,X_r)\,\dd E_r\right |,\\
	&I_3(t)\coloneqq  \left| \int_{t_0}^{t_{\ell_t}}(G(r,X_r)-G(t_{\ell_r},X^\Delta_r))\,\dd B_{E_r} \right|, 
	&&I_4(t)\coloneqq \left| \int_{t_{\ell_t}}^t G(r,X_r)\,\dd B_{E_r} \right|.
\end{align*}
Note that $I_2=I_4=0$ if $t=t_i$ for some $i=0,1,2,...,L$.
Define processes $Z$ and $\widehat Z$ by 
\begin{align}\label{10003}
    Z_t\coloneqq \sup_{t_0\leq s\leq t}|X_s-X_s^\Delta|,\quad\text{and}\quad \widehat Z_t\coloneqq e^{-E_t}Z_t.
\end{align}
That is, $Z$ and $\widehat Z$ are the pathwise uniform error and its exponentially tilted version, respectively.  
Then by \eqref{Eq021-a},
\begin{align*}
    V(t)\coloneqq \E_B\left[Z_t^2\right]\le4\sum_{j=1}^4\E_B\left[\sup_{t_0\le s\le t}|I_j(s)|^2
    \right]
\end{align*}
and
\begin{align}\label{tilt_uniform_error_decomposition}
    \widehat V(t)\coloneqq \E_B\left[\widehat Z_t^2\right]=e^{-2E_t}V(t)\le4\sum_{j=1}^4\E_B\left[e^{-2E_t}\sup_{t_0\le s\le t}\left|I_j(s)\right|^2
    \right].
\end{align}
Next, we examine the right-hand side of \eqref{tilt_uniform_error_decomposition} term by term.

\bigskip
\noindent\textbf{For the $I_1$ term}: For every $r\in[t_0,T]$, conditions \eqref{SDE_condition1-a} and \eqref{SDE_condition4-a} imply
\begin{align}
    |F(r,X_r)-F(t_{\ell_r},X_r^\Delta)|&\le
    |F(r,X_r)-F(t_{\ell_r},X_r)|+|F(t_{\ell_r},X_r)-F(t_{\ell_r},X_r^\Delta)|\notag\\
    &\le K(1+|X_r|)|r-t_{\ell_r}|^\rho +K|X_r-X_r^\Delta|\notag\\
    &\le K Y_T^{(1)}\Delta^\rho+KZ_r,\label{10001}
\end{align}
where $Y_T^{(1)}$ is the one defined in Lemma \ref{lemma_Holder-a}.
Hence, due to the fact $e^{-E_t}\leq 1$,
\begin{align*}
   e^{-E_t}\sup_{t_0\le s\le t}\left| I_1(s)\right|\le e^{-E_t}\int_{t_0}^t\left( K Y_T^{(1)}\Delta^\rho+KZ_r\right)\,\dd E_r\le K Y_T^{(1)}E_T\Delta^\rho + K\int_{t_0}^t e^{-E_t} Z_r\,\dd E_r.
\end{align*}
By the Cauchy--Schwarz inequality and the fact that $\bigl(Y_T^{(1)}\bigr)^2\le 2Y_T^{(2)}$, 
\begin{align*}
    e^{-2E_t}\sup_{t_0\le s\le t}| I_1(s)|^2
    &\le 4K^2E_T^2Y_T^{(2)}\Delta^{2\rho}
    +2K^2\left(\int_{t_0}^t e^{-E_t}Z_r\,\dd E_r\right)^2.
\end{align*}
To deal with the second term on the right hand side, we observe that
\[
\int_{t_0}^t e^{-E_t}Z_r\,\dd E_r
=\int_{t_0}^t e^{-(E_t-E_r)} \widehat Z_r\,\dd E_r
\]
and apply the Cauchy--Schwarz inequality to obtain
\[
e^{-2E_t}\sup_{t_0\le s\le t}|I_1(s)|^2  \le4K^2E_T^2Y_T^{(2)}\Delta^{2\rho}
+ 2K^2\left(\int_{t_0}^te^{-2(E_t-E_r)}\,\dd E_r\right)\left( \int_{t_0}^t\widehat Z^2_r\,\dd E_r\right).
\]
Since $E$ is continuous and nondecreasing, 
\[
\int_{t_0}^te^{-2(E_t-E_r)}\,\dd E_r=\frac1{2}\left(1-e^{-2(E_t-E_{t_0})}\right)\leq \frac1{2},
\]
and hence,
\begin{align}\label{uniform_I1}
    \E_B\left[ e^{-2E_t}\sup_{t_0\le s\le t}|I_1(s)|^2 \right] \le4K^2E_T^2\E_B[Y_T^{(2)}]\Delta^{2\rho}+ K^2\int_{t_0}^t\widehat V(r)\,\dd E_r.
\end{align}

\bigskip
\noindent\textbf{For the $I_2$ term}: By the linear growth condition \eqref{SDE_condition2-a},
\begin{align*}
    \sup_{t_0\le s\le t}|I_2(s)|&\le K Y_T^{(1)} \sup_{t_0\le s\le t} \left(E_s-E_{t_{\ell_s}}\right).
\end{align*}
This implies that 
\begin{align}\label{uniform_I2}
    \E_B\left[e^{-2E_t}\sup_{t_0\le s\le t}|I_2(s)|^2\right] \le 2K^2\E_B[Y_T^{(2)}]\mathcal M_t^2,
\end{align}
where  
\begin{align}\label{10002-0}
    \mathcal M_t
\coloneqq\sup_{t_0\le s\le t} \left(E_s-E_{t_{\ell_s}}\right)= \max\left\{\delta E_0,\ldots,\delta E_{\ell_t-1},R_t\right\},
\end{align}
with 
\begin{align}\label{10002}
    \delta E_i\coloneqq E_{t_{i+1}}-E_{t_i},\quad i=0,1,\ldots,L-1, \quad\text{and}\quad R_t:=E_t-E_{t_{\ell_t}}.
\end{align}

\bigskip
\noindent\textbf{For the $I_3$ term}: 
Recall the Burkholder--Davis--Gundy (BDG) inequality, which states that for any $p>0$, there exists a constant $b_p>0$ such that 
\begin{align}\label{Burkholder}
	\E\left[ \sup_{0\le t\le S} |N_t|^p\right]\le b_p \E\left[ [N,N]_S^{p/2}\right]
\end{align}
for any stopping time $S$ and any continuous local martingale $N$ with quadratic variation $[N,N]$. The constant $b_p$ can be taken independently of $S$ and $N$; see \cite[Chapter 3, Proposition 3.26 \& Theorem 3.28]{KaratzasShreve}.

For every $r\in[t_0,T]$, an argument as in \eqref{10001} yields
\begin{align*}
|G(r,X_r)-G(t_{\ell_r},X_r^\Delta)|
\le K Y_T^{(1)}\Delta^\rho+KZ_r.
\end{align*}
Moreover, observe that the process $N=(N_t)_{t\ge 0}$ defined by $N_t\coloneqq\int_{t_0}^{t_0+t} \left( G(r,X_r)-G(t_{\ell_r},X_r^\Delta) \right)\,\dd B_{E_r}$ is a local martingale, so the BDG inequality \eqref{Burkholder} yields
\begin{align*}
   \E_B\left[ \sup_{t_0\le s\le t}|I_3(s)|^2\right] 
   &=  \E_B\left[ \sup_{t_0\le s\le t}\left|\int_{t_0}^{t_{\ell_s}} \left( G(r,X_r)-G(t_{\ell_r},X_r^\Delta) \right)\,\dd B_{E_r}\right|^2\right]\nonumber\\
    &\le 4\E_B\left[\int_{t_0}^t \left(G(r,X_r)-G(t_{\ell_r},X_r^\Delta)\right)^2\,\dd E_r\right] \nonumber\\
    &\le 8K^2\Delta^{2\rho}\E_B\left[\int_{t_0}^t\bigl(Y_T^{(1)}\bigr)^2\,\dd E_r \right]+8K^2\int_{t_0}^t\E_B[Z_r^2]\,\dd E_r \nonumber\\
    &\le 16K^2E_T\E_B[Y_T^{(2)}]\Delta^{2\rho}+ 8K^2\int_{t_0}^tV(r)\,\dd E_r.
\end{align*}
Therefore, by the obvious inequality $e^{-2E_t}\le e^{-2E_r}$ for any $r\in [t_0,t]$,
\begin{equation}
\label{uniform_I3}
 \E_B\left[e^{-2E_t} \sup_{t_0\le s\le t}|I_3(s)|^2\right] \le 16K^2E_T\E_B[Y_T^{(2)}]\Delta^{2\rho}+ 8K^2\int_{t_0}^t\widehat V(r)\,\dd E_r.
\end{equation}

\bigskip
\noindent\textbf{For the $I_4$ term}: 
Note that whenever $t_i\le s< t_{i+1}$ with $i=0,1,\ldots,L-1$, 
\[
\int_{t_{\ell_s}}^s G(r,X_r)\,\dd B_{E_r}
=\int_{t_i}^s G(r,X_r)\,\dd B_{E_r}
\eqqcolon J_i(s).
\]
{If $t=T$, we set $J_{\ell_T}\equiv 0$.} Then for any fixed integer $q>1$,
\begin{align}
\sup_{t_0\le s\le t}|I_4(s)|^2 
&\le \max\left\{  \max_{0\le i\le \ell_t-1}\sup_{t_i\le s< t_{i+1}}|J_i(s)|^2, \,\sup_{t_{\ell_t}\le s\le t}|J_{\ell_t}(s)|^2\right\}\notag
\\
&=\left(\max\left\{  \max_{0\le i\le \ell_t-1}\sup_{t_i\le s< t_{i+1}}|J_i(s)|^{2q}, \,\sup_{t_{\ell_t}\le s\le t}|J_{\ell_t}(s)|^{2q}\right\}\right)^{1/q}\notag
\\
&\leq\left(\sum_{i=0}^{\ell_t-1}\sup_{t_i\le s< t_{i+1}}|J_i(s)|^{2q}+\sup_{t_{\ell_t}\le s\le t}|J_{\ell_t}(s)|^{2q}\right)^{1/q}.\label{10005}
\end{align}
Since $x\mapsto x^{1/q}$ is concave, we have
\begin{align*}
    \E_B\left[ \sup_{t_0\le s\le t}|I_4(s)|^2 \right]
    \le\Bigg(\sum_{i=0}^{\ell_t-1}\E_B\left[\sup_{t_i\le s< t_{i+1}}|J_i(s)|^{2q}\right]+\E_B\left[ \sup_{t_{\ell_t}\le s\le t}|J_{\ell_t}(s)|^{2q} \right]\Bigg)^{1/q}.
\end{align*}
Applying the BDG inequality \eqref{Burkholder} on each fixed interval $[t_i,t_{i+1})$ and using the linear growth condition \eqref{SDE_condition2-a}, we obtain
\begin{align*}
    \E_B\left[ \sup_{t_i\le s< t_{i+1}}|J_i(s)|^{2q} \right]
    &\le b_q \E_B\left[\left( \int_{t_i}^{t_{i+1}}G(r,X_r)^2\,\dd E_r \right)^q \right]
    \leq b_q\E_B\left[\left( 2K^2Y_T^{(2)}\right)^q\right]\left(\delta E_i\right)^q\\
    &\leq b_q2^{2q-1} K^{2q}\E_B\left[Y_T^{(2q)}\right]\left(\delta E_i\right)^q.
\end{align*}
Similarly, 
\begin{align*}
    \E_B\left[ \sup_{t_{\ell_t}\le s\le t}|J_{\ell_t}(s)|^{2q}\right]&\le b_q2^{2q-1} K^{2q}\E_B\left[Y_T^{(2q)}\right] R_t^q,
\end{align*}
where $R_t$ appears in \eqref{10002}.
Hence, by the obvious inequality $e^{-2E_t}\le 1$, 
\begin{align}\label{uniform_I4}
    \E_B\left[e^{-2E_t} \sup_{t_0\le s\le t}|I_4(s)|^2\right]\le
   C_{K,q} \left(\E_B[Y_T^{(2q)}]\mathcal S_{q,t} \right)^{1/q},
\end{align}
where
\[
    \mathcal S_{q,t}\coloneqq \sum_{i=0}^{\ell_t-1}(\delta E_i)^q+R_t^q
\]
and $C_{K,q}>0$ is some deterministic constant.

\bigskip
\noindent\textbf{Putting together the four estimates}: 
Since $E$ is nondecreasing, both
$t\mapsto \mathcal M_t$ and $t\mapsto \mathcal S_{q,t}$ are nonnegative
and nondecreasing. Thus, with the fact that $R_T=0$, for all $t\in [t_0,T]$,
\begin{equation}
    \label{MS}
   \mathcal M_t\le \mathcal M_T=\max_{0\le i\le L-1}\delta E_i\quad \text{and}\quad
    \mathcal S_{q,t}\le \mathcal S_{q,T}=\sum_{i=0}^{L-1}(\delta E_i)^q.
\end{equation}
Combining \eqref{tilt_uniform_error_decomposition},
\eqref{uniform_I1}, \eqref{uniform_I2}, \eqref{uniform_I3}, \eqref{uniform_I4}, and \eqref{MS},
we obtain
\begin{align}\label{10004}
    \widehat V(t)\le \kappa_q(T)+36K^2\int_{t_0}^t\widehat V(r)\,\dd E_r,
\end{align}
where
\begin{align*}
    \kappa_q(T)&\coloneqq C_K\E_B[Y_T^{(2)}]\Big((E_T^2+E_T)\Delta^{2\rho}+\mathcal M_T^2\Big)+ C_{K,q}
    \left(\E_B[Y_T^{(2q)}]\right)^{1/q}\mathcal S_{q,T}^{1/q}
\end{align*}
for some constants $C_K,C_{K,q}>0$. 
The standard Gronwall inequality for Lebesgue--Stieltjes integrals yields $\widehat V(T)\le\exp\{36K^2E_T \}\kappa_q(T)$, and hence, 
\begin{align}
\label{VT_bound_random_kappa}
    V(T)\le\exp\left\{(36K^2+2)E_T \right\}\kappa_q(T).
\end{align}
Notice that the latter inequality is valid for each sample path of the inverse subordinator, and we now take the full expectation to obtain an upper bound for $\E\left[Z_T^2\right]$. Define the random variable
\begin{align*}
H_{q,T}&:=\exp\{(36K^2+2)E_T\}\Bigg[C_K(1+E_T+E_T^2)\E_B[Y_T^{(2)}]+C_{K,q}
\left(\E_B[Y_T^{(2q)}]\right)^{1/q}\Bigg].
\end{align*}
Then \eqref{VT_bound_random_kappa} implies
\begin{align*}
    \E\left[Z_T^2\right]=\E[V(T)] \le \Delta^{2\rho}\E[H_{q,T}]+\E\left[H_{q,T}\mathcal M_T^2\right]+\E\left[H_{q,T}\mathcal S_{q,T}^{1/q} \right].
\end{align*}
Let $p\coloneqq \frac{q}{q-1}$, so that $\frac 1p+\frac 1q=1$. By H\"older's inequality, the obvious inequality $\max_{1\leq j\leq n} a_j\leq \sum_{j=1}^n a_j$ for $a_j\ge 0$, and the representations in \eqref{MS}, 
\begin{align}
\label{EZ^2}
    \E\left[Z_T^2\right] &\le\Delta^{2\rho}\E[H_{q,T}]+\left( \E[H_{q,T}^p]\right)^{1/p}\left( \E[\mathcal M_T^{2q}] \right)^{1/q}+\left(  \E[H_{q,T}^p] \right)^{1/p} \left( \E\left[  \mathcal S_{q,T}\right]\right)^{1/q}\nonumber\\
    &\le\Delta^{2\rho}\E[H_{q,T}]+\left( \E[H_{q,T}^p] \right)^{1/p}\left(\sum_{i=0}^{L-1} \E[(\delta E_i)^{2q}] \right)^{1/q}+\left(  \E[H_{q,T}^p] \right)^{1/p}\left(\sum_{i=0}^{L-1}\E[(\delta E_i)^q]\right)^{1/q}.
\end{align}

To deal with the terms involving $H_{q,T}$, first observe that 
\[
    (1+E_T+E_T^2)^p\le C_p e^{E_T}
\] 
and 
\[
\left(\E_B[Y_T^{(2)}]\right)^p\le\E_B\left[(Y_T^{(2)})^p\right]\le C_p\E_B[Y_T^{(2p)}], 
\]
where we used Jensen's inequality.
Moreover, since $0<p/q\le 1$ as $q>1$ is an integer, the inequality $x^a\le 1+x$ valid for any $0<a\le 1$ and $x\ge 0$ yields
\[
    \left(\E_B[Y_T^{(2q)}]\right)^{p/q}
    \le
    1+\E_B[Y_T^{(2q)}].
\]
Consequently, Lemma \ref{lemma_Holder-a} and the inequality in \eqref{A.3} 
together imply that
\begin{equation}
\label{HqT}
    \E[H_{q,T}] \le 
    \left(\E[H_{q,T}^p]\right)^{1/p}
    \le C_{\phi,K,q,x_{t_0}}e^{C_{\phi,K,q}T}.
\end{equation}

On the other hand, in terms of the terms involving $\delta E_i$ in \eqref{EZ^2}, 
it follows from \eqref{E-increment-general} that
\begin{align}\label{telescoping_sum_moments}
    \sum_{i=0}^{L-1}
    \E[(\delta E_i)^n]
    \le n! U(\Delta)^{n-1}\sum_{i=0}^{L-1}\left(U(t_{i+1})-U(t_i)\right)
    = n! U(\Delta)^{n-1}(U(T)-U(t_0)).
\end{align}
Applying \eqref{telescoping_sum_moments} with $n=2q$ and $n=q$ in \eqref{EZ^2}, together with \eqref{HqT}, we obtain
\begin{align}\label{uniform_error_general_phi}
  \E\left[Z_T^2\right]  \le
   C_{\phi,K,q,x_{t_0}}e^{C_{\phi,K,q}T}
    \Bigg(\Delta^{2\rho}+\left(U(T)-U(t_0)\right)^{\frac{1}{q}}\left(U(\Delta)^{\frac{2q-1}{q}}+U(\Delta)^{\frac{q-1}{q}} \right)\Bigg).
\end{align}
Since the renewal measure $U(\dd t)$ admits the Laplace transform $\lambda\to \phi(\lambda)^{-1}$ as observed in \eqref{U-Laplace}, it follows that 
\begin{align}\label{U-exp-bnd}
     \frac{1}{\phi(\lambda)}
    =\int_0^\infty e^{-\lambda t}U(dt)
    \ge \int_{t_0}^T e^{-\lambda t}U(dt)
    \ge e^{-\lambda T}(U(T)-U(t_0)),\quad \lambda>0,
\end{align}
and therefore, $\left(U(T)-U(t_0)\right)^{\frac{1}{q}}\le C_{q,\phi}e^{C_{q}T}$.
Using this bound and \eqref{U-and-phi} in \eqref{uniform_error_general_phi} yields 
\begin{align}\label{1511}
\E\left[\sup_{t_0\le t\le T}|X_t-X_t^\Delta|^2\right] \le  C_1e^{C_2T}
   \Bigg(\Delta^{2\rho}+\frac{1}{\phi(1/\Delta)^{1-1/q}} \Bigg),\quad \Delta\in(0,1).
\end{align}

Finally, fix $\varepsilon\in(0,1)$ and choose an integer $q>1$ such that $1/q\le\varepsilon$. Since $\phi$ is increasing and $\Delta\in(0,1)$, it follows that $\phi(1/\Delta)\ge\phi(1)>0$. Hence,
\[
\frac{1}{\phi(1/\Delta)^{1-1/q}}=\phi(1/\Delta)^{1/q-\varepsilon}\frac{1}{\phi(1/\Delta)^{1-\varepsilon}}\le\phi(1)^{1/q-\varepsilon}\frac{1}{\phi(1/\Delta)^{1-\varepsilon}}.
\]
Absorbing the factor $\phi(1)^{1/q-\varepsilon}$ into the constant in \eqref{1511}
gives the desired statement in \eqref{1427}.

If $\phi$ additionally satisfies the property that $\phi(\lambda)\ge c\lambda^\alpha$ for
$\lambda\ge\lambda_0$ with $\alpha\in(0,1]$, then, by retaking the constant 
$c>0$ if necessary, we may assume that the same inequality holds
for all $\lambda\ge1$. Therefore,
\[
    \frac{1}{\phi(1/\Delta)^{1-\varepsilon}}
    \le
    C\Delta^{\alpha(1-\varepsilon)},
    \qquad \Delta\in(0,1),
\]
which, together with \eqref{1427}, yields \eqref{1509}.
\end{proof}

\begin{corollary}
Under the assumptions of Theorem \ref{thm:strong-approx} with $G\equiv 0$, for every $\varepsilon\in(0,1)$, there exists a constant $C_1=C_1(x_{t_0},K,\phi,\varepsilon)>0$, independent of $\Delta$
and $T$, such that
\begin{align*}
\E\left[\sup_{t_0\le t\le T}|X_t-X_t^\Delta|^2\right]\le C_1e^{C_1T}\left(\Delta^{2\rho}+\frac{1}{\phi(1/\Delta)^{2-\varepsilon}}\right),\qquad \Delta\in(0,1).
\end{align*}
Moreover, if there exist constants $\alpha\in(0,1]$ and $c,\lambda_0>0$
such that $\phi(\lambda)\ge c\lambda^\alpha$ for all $\lambda\ge\lambda_0$, then there exists a constant $C_2=C_2(x_{t_0},K,\phi,\varepsilon)>0$, independent of $\Delta$
and $T$, such that
\[
\E\left[\sup_{t_0\le t\le T}|X_t-X_t^\Delta|^2\right]\le C_2e^{C_2T}\Delta^{\min\{2\rho,\alpha(2-\varepsilon)\}},\qquad \Delta\in(0,1).
\]
 Consequently, $X^\Delta$ converges strongly to $X$ uniformly over
$[t_0,T]$ with order $\eta=\min\left\{\rho,\alpha\left(1-\frac{\varepsilon}{2}\right)\right\}$.
\end{corollary}

\begin{proof}
The first statement follows from the proof of Theorem \ref{thm:strong-approx} by setting $I_3(t)=I_4(t)=0$, which makes the last term of \eqref{EZ^2} disappear. 
\end{proof}

\begin{remark}\label{r:rate0}
\begin{enumerate}[(a)]
\item If the subordinator in Theorem \ref{thm:strong-approx} has positive drift, then representation \eqref{bernstein-exp} immediately gives $\phi(\lambda)\ge \mathfrak{b}\lambda$ for all $\lambda>0$, and consequently,  \eqref{1509} holds with $\alpha=1$. 

\item The constant $C_1$ in \eqref{1427} generally diverges to $\infty$ as $\varepsilon\downarrow0$. Indeed, in the proof of Theorem \ref{thm:strong-approx}, the integer $q$ is chosen so that $1/q\le\varepsilon$, so $q\to\infty$ as $\varepsilon\downarrow0$, whereas the constant $C_{\phi,K,q,x_{t_0}}$ in \eqref{uniform_error_general_phi}
contains the factor $((2q)!)^{1/q}$, which diverges as $q\to\infty$. Therefore, one cannot simply let $\varepsilon\downarrow0$ in \eqref{1427} to obtain an upper bound of the form
$Ce^{CT}\left(\Delta^{2\rho}+\phi(1/\Delta)^{-1}\right)$.
In fact, in Section \ref{s:sharp}, we will observe that this is not due to a technical issue of the proof of Theorem \ref{thm:strong-approx} but due to an intrinsic property of the piecewise-constant interpolation. More precisely, we will provide an example for which the $L^2$ convergence rate cannot be improved from $\Delta^{2\rho}+\phi(1/\Delta)^{\varepsilon-1}$ to $\Delta^{2\rho}+\phi(1/\Delta)^{-1}$.
\end{enumerate}
\end{remark}

\begin{remark}\label{rem:Zuo}
A statement similar in nature to Theorem \ref{thm:strong-approx} has appeared as \cite[Theorem 3.1]{Zuo2026}; however, the approach taken and the scope of the statement given in \cite[Theorem 3.1]{Zuo2026} are clearly different from those in Theorem \ref{thm:strong-approx} of the present paper. 

In fact, \cite[Theorem 3.1]{Zuo2026} is restricted to the setting where the time change is an inverse $\alpha$-\textit{stable} subordinator with $\alpha\in(\frac 12,1)$ and the SDE coefficients are homogeneous (i.e., they depend on $X_t$ only). The argument used in the proof calls for a statement on the modulus of continuity for inverse \textit{stable} subordinators appearing in \cite{IKMS17}, which allows the quantity
\[
\mathcal M_t
=\sup_{t_0\le s\le t} \left(E_s-E_{t_{\ell_s}}\right) 
\]
to be well-controlled. 
In contrast, our statement, Theorem \ref{thm:strong-approx}, is applicable to a much broader class of inverse subordinators, including the inverses of $\alpha$-stable subordinators with $\alpha\in(0,1)$ and non-stable subordinators. 
We have overcome the hurdle for this broad generalization by (i) the exponential tilting in \eqref{10003}, and (ii) the simple estimate
\begin{align}\label{eq:simple estimate}
    \mathcal M_T^{2q}
    =\max_{0\le i\le L-1}(\delta E_i)^{2q}
    \le \sum_{i=0}^{L-1}(\delta E_i)^{2q}
\end{align}
and a similar estimate in \eqref{10005}. 

For (i), the reader should note that, without the exponential tilting of $Z$, one could obtain an inequality similar to \eqref{10004} but \textit{with the coefficient of the integral taking the form $f(E_T)$ for some polynomial $f$ of degree 2,} as in the proof of \cite[Theorem 3.1]{Zuo2026}. Using the Gronwall inequality would then yield a factor $e^{Cf(E_t)}\ge e^{c E_t^2}$, which would restrict the validity of the final estimate to a limited range of $\alpha$ values due to Theorem \ref{theorem_expo_moment}.
Regarding (ii), let us emphasize that our seemingly crude estimates do not contribute to lowering the order of strong convergence derived in \cite[Theorem 3.1]{Zuo2026} for a limited class of inverse stable subordinators; in fact, the order remains the same. 
However, a much more significant benefit of using the simple estimate \eqref{eq:simple estimate} in \eqref{10005} is that it enables us to fully avoid the discussion of the modulus of continuity associated with the stochastic integral $M_t\coloneqq\int_{t_0}^t G(D_{r-},X_{D_{r-}})\, dB_r$, which is not a trivial question. In fact, with the notations used in the present paper and with $t_0=0$ for simplicity, the proof of \cite[Theorem 3.1]{Zuo2026} claims that the inequality 
\[
    \sup_{0\le s\le t}|I_4(s)|
    \le \sup_{0\le r<u\le E_T,\; u-r \le \Xi_\Delta} |M_u-M_r|
\]
together with \cite[Theorem 7]{JinKobayashi2019} implies 
\[
    \E_B\left[\sup_{0\le s\le t}|I_4(s)|^2\right]
    \le C_p \E_B[Y^{(2)}_T] \Xi_\Delta \log\left(\frac{4E_T}{\Xi_\Delta}\right),
\]
where $\Xi_\Delta$ denotes the maximal increment of $E$ over $[0,T]$ defined by 
\[
    \Xi_\Delta\coloneqq \sup_{0\le u\le T-\Delta} (E_{u+\Delta}-E_u);
\]
however, this implication is not justified there. 
In our view, this is not an immediate consequence of \cite[Theorem 1]{FischerNappo}, which plays an important role in the proof of \cite[Theorem 7]{JinKobayashi2019}. 

\end{remark}

\section{Strong approximation of time-changed SDEs with piecewise-constant interpolation: the \textit{synchronous} case}\label{s:approx-2}

The goal of this section is to study the convergence order of a strong approximation scheme for time-changed SDEs of the form 
\begin{align}\label{duality-X}
    \dd X_t=F(E_t,X_t)\dd E_t+G(E_t,X_t) \dd B_{E_t},\quad t\ge0,
\end{align}
with initial condition $X_0=x_0$, where the coefficients depend on the inverse subordinator $E_t$ directly rather than on the time variable $t$. We say the coefficients are \textit{synchronous} because $E_t$ appearing in the coefficients shares the same constant periods with the driving processes $E_t$ and $B_{E_t}$.

The so-called \textit{duality principle} established in \cite[Theorem 4.2]{Kobayashi2011} has played an important role in various analyses of time-changed SDEs with synchronous coefficients. Namely, if the process $X$ satisfies SDE \eqref{duality-X}, then the process $Y$ defined by $Y_t=X_{D_t}$, where $D$ is the underlying subordinator, satisfies the classical It\^o SDE 
\begin{align}\label{duality-Z}
    \dd Y_t=F(t,Y_t)\dd t+G(t,Y_t)\dd B_t, \quad t\ge0,
\end{align}
with initial condition $Y_0=x_0$. Similarly, if $Y$ satisfies the SDE \eqref{duality-Z}, then the process $X$ defined by $X_t=Y_{E_t}$  satisfies SDE \eqref{duality-X}. 
The duality principle is particularly useful when the underlying Brownian motion $B$ and the inverse subordinator $E$ are independent, as it makes it possible to analyze the time-changed SDE in \eqref{duality-X} by separately analyzing the inverse subordinator $E$ and the classical It\^o SDE in \eqref{duality-Z}. 
In this section, we will utilize the duality principle to analyze the strong convergence order for SDE \eqref{duality-X}.

The class of time-changed SDEs with synchronous coefficients contains the important subclass with \textit{homogeneous} coefficients; i.e., $F(E_t,X_t)=F(X_t)$ and $G(E_t,X_t)=G(X_t)$.
In fact, this can also be considered a subclass of time-changed SDEs of the form \eqref{SDE_001-a} with asynchronous coefficients already discussed in Section \ref{s:approx-1}. 
Therefore, for numerical approximation of time-changed SDEs in this particular subclass, one can take two different approaches in the analysis of convergence rates: one approach taken in Section \ref{s:approx-1} that utilized exponential tilting, and the other to be taken in the present section based on the duality principle.  
In short, for this subclass of time-changed SDEs, our analysis in this section shows no improvement in the convergence rate from what we observed in Section \ref{s:approx-1}. Let us also emphasize that, as mentioned in Remark \ref{r:rate0}, the obtained convergence rate is sharp and cannot be improved; we defer this important observation to Section \ref{s:sharp}. 

Let us now clarify the setting of this section. We study time-changed SDEs of the form 
\begin{align}\label{SDE_002}
    X_t=x_{0}+\int_{0}^tF(E_r,X_r)\dd E_r+\int_{0}^tG(E_r,X_r)\dd B_{E_r},\quad t\ge0,
\end{align}
through the corresponding It\^o SDEs of the form 
\begin{align}\label{SDE_002-Y}
    Y_t=y_{0}+\int_{0}^tF(r,Y_r)\dd r+\int_{0}^tG(r,Y_r)\dd B_r,\quad t\ge0,
\end{align}
where the functions $F(t,x)$ and $G(t,x)$ satisfy the global version of the conditions \eqref{SDE_condition1-a}--\eqref{SDE_condition4-a} with $t_0=0$; that is, the coefficients $F(t,x)$ and $G(t,x)$ satisfy conditions \eqref{SDE_condition1-a}--\eqref{SDE_condition4-a} uniformly in $[0,\infty)\times \R^d$. 
For the proofs in this section, we keep dimensions $m=d=1$ 
as in Section \ref{s:approx-1} for simplicity of discussions and expressions. 
Additionally, we assume that the Laplace exponent $\phi$ of the subordinator $D$ satisfies $\phi(\infty)=\infty$ so that $D$ has strictly increasing paths, and hence,  its inverse $E$ has continuous, nondecreasing paths. {Furthermore, the Brownian motion $B$ and the inverse subordinator $E$ are assumed to be independent, which in turn makes the solution $Y$ to the classical It\^o SDE \eqref{SDE_002-Y} independent of $E$.}

Note that we take the initial time of the SDE \eqref{SDE_002} to be $t_0=0$ so that the duality principle can be applied. 
In fact, when $t_0>0$, the duality principle as described above no longer holds, and even though a modification of the duality principle is possible, the resulting SDE will require an additional layer of careful analysis. 
We clarify this point in Remark \ref{r:t_0>0}. 

With the duality principle and the synchronization approach presented in \cite{Kobayashi2011} in mind, throughout this section, we fix the representation $X_t=Y_{E_t}$ for the solution $X$ to SDE \eqref{SDE_002}, with $Y$ being the corresponding classical It\^o process in \eqref{SDE_002-Y} with $y_0=x_0$, and study the Euler--Maruyama-type approximation of $X$ through the classical Euler--Maruyama approximation of $Y$. 
As in Section \ref{s:approx-1}, $X^\Delta=(X^\Delta_t)_{t\in[0,T]}$ denotes the approximation process defined {in the same spirit as in \eqref{Approx_001-a}--\eqref{Approx_003-a}. Namely, fix a time-grid as in \eqref{deterministic-time-a} with $t_0=0$, set
\begin{align}
	&X^\Delta_{0}
	=x_{0}\label{Approx_001-b},\\
	&X^\Delta_{t_{\ell+1}}\!\!
=X^\Delta_{t_{\ell}}+F(E_{t_{\ell}},X^\Delta_{t_{\ell}})(E_{t_{\ell+1}}-E_{t_{\ell}})+G(E_{t_{\ell}},X^\Delta_{t_{\ell}})(B_{E_{t_{\ell+1}}}-B_{E_{t_{\ell}}})
	\label{Approx_002-b}
\end{align}
for $\ell=0,1,2,\ldots,L-1$, and use piecewise-constant interpolation as in \eqref{Approx_003-a} to define the process $X^\Delta$ in continuous time, with $\ell_t$ defined in \eqref{1505}.}

We now state the main theorem of this section. The statement itself looks similar to the statement of Theorem \ref{thm:strong-approx}, but the approach to be taken is essentially different due to the synchronous nature of the coefficients with respect to the driving processes.

\begin{theorem}\label{thm:synchro-strong}
Let $E$ be the inverse of a subordinator $D$ with Laplace exponent $\phi$ in \eqref{bernstein-exp} with $\mathfrak{b}>0$ or $\nu(0,\infty)=\infty$. 
Let $B$ be a Brownian motion independent of $D$.
Let $X$ be the unique solution of SDE \eqref{SDE_002} on $[0,T]$, where the coefficients $F(t,x)$ and $G(t,x)$ satisfy conditions \eqref{SDE_condition1-a}--\eqref{SDE_condition4-a} uniformly in $[0,\infty)\times \R^d$. 
Let $X^\Delta$ be the approximation process constructed in \eqref{Approx_001-b}, \eqref{Approx_002-b} and \eqref{Approx_003-a}
based on deterministic time points \eqref{deterministic-time-a} with $t_0=0$ and $L\ge 2$.
Then, for every $\varepsilon\in (0,\min\{1,2\rho\})$, there exists a constant $C_1=C_1(\phi,\varepsilon,K,\rho,x_0)>0$, independent of $T$ and $\Delta$, such that
    \begin{align}\label{1607-new}
        \E\left[\sup_{0\le s\le T}|X_s- X^\Delta_s|^2\right]&\le 
        {C_1e^{C_1T}\left(\frac{1}{\phi(1/\Delta)}\right)^{\min\{1,2\rho\}-\varepsilon}, 
        \quad \Delta\in(0,1).}
    \end{align}
Moreover, if there exist constants $\alpha\in(0,1]$ and $c,\lambda_0>0$ such that $\phi(\lambda)\ge c\lambda^\alpha$ for all $\lambda\ge\lambda_0$, then there exists a constant $C_2=C_2(\phi,\varepsilon,K,\rho,x_0)>0$, independent of $\Delta$
and $T$, such that
        \begin{align*}
        \E\left[\sup_{0\le s\le T}|X_s- X^\Delta_s|^2\right]&\le 
        {C_2e^{C_2T}\Delta^{\alpha(\min\{1,2\rho\}-\varepsilon)}, 
        \quad \Delta\in(0,1).}
    \end{align*}
  Consequently, $X^\Delta$ converges strongly to $X$ uniformly over $[0,T]$ with order $\eta=\frac{\alpha}{2}(\min\{1,2\rho\}-\varepsilon)$.

\end{theorem}

The proof of Theorem \ref{thm:synchro-strong} requires some preparations that enable us to exploit the duality principle and the moduli of continuity of the SDEs. This approach makes the proof of Theorem \ref{thm:synchro-strong} neat and elegant.

First, note that due to the representation $X_t=Y_{E_t}$, the deterministic time horizon $T$ for $X$ corresponds to the random time $E_T$ for $Y$, and the deterministic time step $t_n-t_{n-1}$ for $X$ corresponds to the random time step $E_{t_n}-E_{t_{n-1}}$ for $Y$. Since $E_t$ for any fixed $t>0$ is supported in $(0,\infty)$, the latter motivates us to study the process $Y$ over the infinite time interval $[0,\infty)$. 
Now, for the It\^o process $Y$, consider an arbitrary mesh
\begin{align}\label{1503}
    {0=s_0\le s_1\le \cdots\le s_n\nearrow\infty,}
\end{align}
and let
\begin{align}\label{1502}
	n_{\{s_n\},s}\coloneqq n_s\coloneqq \max\{n\in \mathbb{N}\cup \{0\}: s_n\le s\}.
\end{align}
Let $\wt Y_0\coloneqq y_0$ and
\begin{align}\label{EM-Y}
    \wt Y_s\coloneqq \wt Y_{s_n}+(s-s_n)F(s_n,\wt Y_{s_n})+G(s_n,\wt Y_{s_n})(B_s-B_{s_n}),\quad s\in (s_n,s_{n+1}].
\end{align}
While analyzing the strong error of the approximation of $Y$, it is essential to keep track of all dependencies on (i) the time horizon $T$ and (ii) the maximum size of the time steps given by 
\begin{align}\label{1504}
    \Delta_{\{s_n\}}(s)\coloneqq \max\left\{\max\left\{s_{n+1}-s_{n}:n\le n_{s}-1\right\},s-s_{n_{s}}\right\},
\end{align}
which includes the size of the last step from $s_{n_{s}}$ to $s$, with $n_s$ defined in \eqref{1502}.

The following proposition expresses an upper bound in terms of $\Delta_{\{s_n\}}(s)$ for the $L^2$ strong error ascribed to the use of the approximation process $\wt Y$ in \eqref{EM-Y}. 

\begin{proposition}\label{JumKobayashi}
Let $Y$ be the unique solution of SDE \eqref{SDE_002-Y}, where ${y_0\in \R^d}$ 
and the coefficients $F(t,x)$ and $G(t,x)$ satisfy {the conditions \eqref{SDE_condition1-a}--\eqref{SDE_condition4-a} uniformly in $[0,\infty)\times \R^d$.} Let $\wt Y$ be the approximation process constructed in \eqref{EM-Y} based on deterministic time points \eqref{1503}. 
Then, there exists a constant $C>0$ 
such that for all $s\ge0$,
\begin{align}\label{1501}
    \E\left[\sup_{0\le r\le s}|Y_r-\wt Y_r|^2\right]\le C(1+|y_0|^2)e^{Cs}\max\{\Delta_{\{s_n\}}(s),\Delta_{\{s_n\}}(s)^{2\rho}\}.
\end{align}
\end{proposition}

\begin{proof}
This statement is an improvement of \cite[Proposition 3.1]{JumKobayashi} and follows from the proof of \cite[Lemma 5.1]{BCT2026}, where an upper bound of the form $Ce^{Cs}h^{\min\{1,2\rho\}}$ is obtained, with $h\coloneqq \sup\{s_n-s_{n-1}:n\in\N\}\in(0,1)$.
However, the proof can be rewritten in terms of the maximum step size $\Delta_{\{s_n\}}(s)$, which results in \eqref{1501}.
\end{proof}

We will also need a lemma dealing with maximum fluctuations of It\^o diffusions and their time-changed versions.
\begin{lemma}\label{1314}
    Assume that $F(t,x)$ and $G(t,x)$ satisfy {the conditions \eqref{SDE_condition1-a}--\eqref{SDE_condition4-a} uniformly in $[0,\infty)\times \R^d$.} Assume further that 
    $Y$ and $X$ are the strong solutions of SDEs \eqref{SDE_002-Y} and \eqref{SDE_002} with {$y_0\in\R^d$ and $x_0\in\R^d$}, respectively.
    Then  for any $p\in[1,\infty)$, there exist constants $C_1,C_2>0$  
    such that for all $s\ge0$ and $t\ge 0$, 
    \begin{align}
\E\left[1+\sup_{r\in[0,s]}|Y_r|^p\right]\le C_1(1+|y_0|^p)e^{ C_1s}\label{1314a}
        \end{align}
        and
        \begin{align}
\E\left[1+\sup_{r\in[0,t]}|X_r|^p\right]\le C_2(1+|x_0|^p)e^{C_2t}.\label{1314b}
    \end{align}
\end{lemma}
\begin{proof}
    We may regard \eqref{1314a} as the special case of Lemma \ref{lemma_Holder-a} with $D_t=t$, or equivalently, $E_t=t$.
    Combining the independence between $Y$ and $E$ and the inequalities in \eqref{1314a} and \eqref{A.3},
    \begin{align*}
\E\left[1+\sup_{r\in[0,t]}|X_r|^p\right]
&=\E\left[1+\sup_{r\in[0,E_t]}|Y_r|^p\right]
        \le C_1(1+|x_0|^p)\E[e^{C_1 E_t}]\\
        &\le C_2(1+|x_0|^p)e^{C_2t}.
    \end{align*}
This is \eqref{1314b}, and the proof is complete.
\end{proof}

Recall that our goal is to analyze the strong error of the use of the approximation $X^\Delta$ in \eqref{Approx_001-b}--\eqref{Approx_002-b}, with piecewise-constant interpolation in \eqref{Approx_003-a}, in approximating the solution $X$ to SDE \eqref{SDE_002}. However, for the purpose of the analysis, we will employ the approximation $\wt X=(\wt X_t)_{t\in[0,T]}$ \textit{with continuous interpolation}; namely, $\wt X_0=x_0$, and 
\begin{align}\label{1512}
    \wt X_t=\wt X_{t_\ell}+(E_t-E_{t_\ell})F(E_{t_\ell},\wt X_{t_\ell})+G(E_{t_\ell},\wt X_{t_\ell})(B_{E_t}-B_{E_{t_\ell}}),\quad t\in (t_\ell,t_{\ell+1}].
\end{align}
Note that the sample paths of $\wt X$ coincide with those of the time-changed process $(\wt Y_{E_t})_{t\ge 0}$ if the mesh for $\wt Y$ in \eqref{1503} is chosen as $s_n\coloneqq E_{t_n}$. On the other hand, $\wt X_{t_n}=X^\Delta_{t_n}$ at each time-step point since the two approximation processes differ only in the interpolation methods.

Establishing the following statement, which concerns the strong approximation of $X$ by $\wt X$, is a crucial step in analyzing the error induced by the use of $X^\Delta$. 
Proposition \ref{JumKobayashi} together with the duality principle in \cite{Kobayashi2011} allows us to prove this important statement. 

\begin{theorem}\label{thm:approx-synchr-conti}
    Under the assumptions of Theorem \ref{thm:synchro-strong}, let $\wt X$ be the approximation process defined in \eqref{1512}. 
    Then, for every $\varepsilon\in (0,\min\{1,2\rho\})$, there exists a constant $C_1=C_1(\phi,\varepsilon,K,\rho,x_0)>0$, independent of $\Delta$ and $T$, such that
    \begin{align}\label{1607}
        \E\left[\sup_{0\le s\le T}|X_s-\wt X_s|^2\right]\le 
        C_1e^{C_1T}\left(\frac{1}{\phi(1/\Delta)}\right)^{\min\{1,2\rho\}-\varepsilon},
        \quad \Delta\in(0,1).
    \end{align}
    Moreover, if there exist constants $\alpha\in(0,1]$ and $c,\lambda_0>0$ such that $\phi(\lambda)\ge c\lambda^\alpha$ for all $\lambda\ge\lambda_0$, then there exists a constant $C_2=C_2(\phi,\varepsilon,K,\rho,x_0)>0$, independent of $\Delta$ and $T$, such that
        \begin{align}\label{1517}
        \E\left[\sup_{0\le s\le T}|X_s- \wt X_s|^2\right]&\le 
        C_2e^{C_2T}\Delta^{\alpha(\min\{1,2\rho\}-\varepsilon)}, 
        \quad \Delta\in(0,1).
    \end{align}
Consequently, $\wt X$ converges strongly to $X$ uniformly over $[0,T]$ with order $\eta=\frac{\alpha}{2}(\min\{1,2\rho\}-\varepsilon)$.
\end{theorem}

\begin{proof}
{In the proof, for any parameter $\theta$, the notation $C_\theta$ represents a positive constant depending on $\theta$, which may vary from line to line.}

    Since $B$ is assumed to be independent of $E$, so are both $Y$ and $\wt Y$, where $\wt Y$ is defined in \eqref{EM-Y}. 
   Take the mesh for $\wt Y$ in \eqref{1503} as $s_n\coloneqq E_{t_n}$ so that $\wt X_t=\wt Y_{E_t}$ for all $t$. Then by Proposition \ref{JumKobayashi},
    \begin{align*}
        \E\left[\sup_{0\le s\le T}|X_s-\wt X_s|^2\right]
        &=\E\left[\sup_{0\le s\le T}|Y_{E_s}-\wt Y_{E_s}|^2\right]=\E\left[\sup_{0\le s\le E_T}|Y_{s}-\wt Y_{s}|^2\right]\\
        &\le  \E\left[C_{K,\rho}(1+|x_0|^2)e^{C_{K,\rho}E_T} \max\{\Delta_{\{E_{t_n}\}}(E_T),\Delta_{\{E_{t_n}\}}(E_T)^{2\rho}\}\right].
    \end{align*}
Observe that by \eqref{1502} and \eqref{1505}, 
    \[
        n_{\{E_{t_n}\},E_T}
        = \max\{n\in \mathbb{N}\cup \{0\}: E_{t_n}\le E_T\}
        =\max\{n\in \mathbb{N}\cup \{0\}: t_n\le T\}
        =\ell_T=L,
    \]
    which, together with \eqref{1504}, gives
    \[
        \Delta_{\{E_{t_n}\}}(E_T)  
        =\max\left\{\delta E_n:n\le n_{\{E_{t_n}\},{E_T}}-1\right\}
        =\max\left\{\delta E_n:n\le L-1\right\}
        =\MM_T,
    \]
    with $\delta E_i= E_{t_{i+1}}-E_{t_i}$ and $\MM_T$ defined in \eqref{10002} and \eqref{10002-0}, respectively. 
    Thus, we have obtained 
    \begin{align}   
    \E\left[\sup_{0\le s\le T}|X_s-\wt X_s|^2\right]
    &\le\E\left[C_{K,\rho}(1+|x_0|^2)e^{C_{K,\rho}E_T}\max\{\MM_T,\MM_T^{2\rho}\}\right]\notag\\
    &\le C(1+|x_0|^2)\E\left[e^{pC_{K,\rho}E_T}\right]^{1/p}\E\left[\max\{\MM_T^q,\MM_T^{2\rho q}\}\right]^{1/q},\label{1404} 
    \end{align}    
    where we applied H\"older's inequality with $p,q>1$ and $1/p+1/q=1$. (The constant $q$ will be specified later.)
    By the inequality in \eqref{A.3}, 
    the factor $\E\left[e^{pC_{K,\rho}E_T}\right]^{1/p}$ is finite and bounded by an expression of the form $C_pe^{C_p T}$. For the second factor in \eqref{1404}, 
    \begin{align}
        \E\left[\max\{\MM_T^q,\MM_T^{2\rho q}\}\right]^{1/q}
        &\le \left(\sum_{i=0}^{L-1}(\E[(\delta E_i)^q]+\E[(\delta E_i)^{2\rho q}])\right)^{1/q}\notag\\
        &\le \left(\sum_{i=0}^{L-1}\E[(\delta E_i)^q]\right)^{1/q}+\left(\sum_{i=0}^{L-1}\E[(\delta E_i)^{2\rho q}]\right)^{1/q}.\label{1507}
    \end{align}

    {Now, for a fixed $\varepsilon\in(0,\min\{1,2\rho\})$, choose an integer $q>\max\{1/\varepsilon, 1/(2\rho)\}$. 
    By Theorem \ref{Theorem_increments}(b)(i),
    \begin{align}
        \sum_{i=0}^{L-1}\E[(\delta E_i)^q]
        &\le C_q  \sum_{i=0}^{L-1}\big(U(t_{i+1})-U(t_i)\big)U(t_{i+1}-t_i)^{q-1}\notag\\
        &\le C_q U(T)U(\Delta)^{q-1}
        =C_q \left(\frac{U(T)}{U(\Delta)}\right)U(\Delta)^q,\label{1507+}
    \end{align}
    where the term $U(\Delta)$ is obtained by recalling that $t_{i+1}-t_i\le \Delta$ for $i=0,1,\ldots,L-1$ and noticing that $t\mapsto U(t)$ is nondecreasing. Since $q>1/\varepsilon$ and $U(\Delta)\le U(T)$ for $\Delta\in(0,1)$,
    \begin{align}\label{1602-a}
         \left(\sum_{i=0}^{L-1}\E[(\delta E_i)^q]\right)^{1/q}
         \le C_q \left(\frac{U(T)}{U(\Delta)}\right)^{1/q}U(\Delta)
         \le C_q \left(\frac{U(T)}{U(\Delta)}\right)^{\varepsilon}U(\Delta).
    \end{align}
To deal with the second term in \eqref{1507}, let
$n=\lceil 2\rho q\rceil$ and $\theta=(2\rho q-1)/(n-1)\in(0,1]$. Then $n\ge 2$ since $q>1/(2\rho)$.
By H\"older's inequality and Theorem \ref{Theorem_increments},
\begin{align*}
\E[(\delta E_i)^{2\rho q}]
&=\E[(\delta E_i)^{1-\theta}(\delta E_i)^{n\theta}]
\le\left(\E[\delta E_i]\right)^{1-\theta}
\left(\E[(\delta E_i)^n]\right)^\theta\\
&\le(n!)^\theta\bigl(U(t_{i+1})-U(t_i)\bigr)U(\Delta)^{(n-1)\theta}\\
&=(n!)^\theta\bigl(U(t_{i+1})-U(t_i)\bigr)U(\Delta)^{2\rho q-1}.
\end{align*}
Consequently,
\begin{align}\label{1602-b}
\left(\sum_{i=0}^{L-1}
\E[(\delta E_i)^{2\rho q}]\right)^{1/q}
\le 
C_{\rho,q}\left(\frac{U(T)}{U(\Delta)}\right)^{1/q}U(\Delta)^{2\rho}
\le C_{\rho,q}\left(\frac{U(T)}{U(\Delta)}\right)^{\varepsilon}U(\Delta)^{2\rho}.
\end{align}
Recall from \eqref{U-exp-bnd} and \eqref{U-and-phi} that $U(T)\le C_\phi e^{T}$ and $U(\Delta)\le e\phi(1/\Delta)^{-1}$, respectively. 
These two bounds, together with \eqref{1404}, \eqref{1507}, \eqref{1602-a} and \eqref{1602-b}, yield \eqref{1607}.}
\end{proof}

With Lemma \ref{1314} and Theorem \ref{thm:approx-synchr-conti} at hand, we are now ready to prove the main theorem of this section, Theorem \ref{thm:synchro-strong}.

\begin{proof}[Proof of Theorem \ref{thm:synchro-strong}]
   In the proof, for any parameter $\theta$, the notation $C_\theta$ represents a positive constant depending on $\theta$, which may vary from line to line.
   
   Recall that for $t\in [t_n,t_{n+1})$, it holds that $X^\Delta_t=X^\Delta_{t_n}=\wt X_{t_n}$ which implies
    \begin{align*}
        |X_t-X^\Delta_t|^{2}=|X_t-\wt X_{t_n}|^{2}\le 2|X_{t_n}-\wt X_{t_n}|^{2}+2|X_t- X_{t_n}|^{2}.
    \end{align*}
    Therefore,
    \begin{align}
        \sup_{0\le s\le T}|X_s-X^\Delta_s|^{2}&\le 2\sup_{0\le s\le T}|X_s-\wt X_s|^{2}+2\max_{i=0,\ldots,L-1}\left\{\sup_{t_{i}\le s < t_{i+1} }|X_{s}-X_{t_i}|^2\right\}\notag\\
        &\eqqcolon I_1(T)+I_2(T).\label{I1+I2}
    \end{align}
    
    For the $I_1(T)$ term,
    by Theorem \ref{thm:approx-synchr-conti}, for all $\varepsilon\in (0,\min\{1,2\rho\})$, there exists $C=C(\phi,\varepsilon,K,\rho,x_0)>0$ such that
    \begin{align}\label{1113-a}
        \E [I_1(T)]
        \le C(1+|x_0|^2)e^{C\,T}\left(\frac{1}{\phi(1/\Delta)}\right)^{\min\{1,2\rho\}-\varepsilon}.
    \end{align}
     For the $I_2(T)$ term, note that 
    \begin{align*}
        I_2(T)
        =2\max_{i=0,\ldots,L-1}\left\{\sup_{t_{i}\le s < t_{i+1} }|Y_{E_s}-Y_{E_{t_i}}|^2\right\}
        \le 2\sup_{r,s\in [0,E_T],|r-s|\le \MM_T}|Y_r-Y_s|^2.
    \end{align*}
{Then a close inspection of the assumptions and proof of \cite[Theorem 1]{FischerNappo} together with Lemma \ref{1314} implies that 
    \begin{align}\label{1506}
        \E[I_2(T)]\le 2\E\left[\sup_{r,s\in [0,E_T],|r-s|\le \MM_T}|Y_r-Y_s|^2\right]\le \E\left[f_Y(E_T)\MM_T\log\frac{2E_T}{\MM_T}\right],
    \end{align}   
    where 
    \begin{align}\label{1506+}
        f_Y(E_T)\coloneqq C_{K,{\gamma}}\E_B\left[\sup_{s\in [0,E_T]}(1+|Y_s|^{2+\gamma})\right]\le C_{K,{\gamma}} (1+|x_0|^{2+\gamma}) e^{\tilde c_K E_T},
    \end{align}
    {with $\gamma>0$.} Notice that $f_Y(E_T)$ involves the $(2+\gamma)$-th {power rather than the square} because the proof of \cite[Theorem 1]{FischerNappo} relies on H\"older's inequality.

Fix $\varepsilon\in(0,\min\{1,2\rho\})$ and choose an
integer $q>2/\varepsilon$. Let
$p=q/(q-1)$. By H\"older's inequality,
\begin{align}
\label{1442}
\E\left[f_Y(E_T)\MM_T\log\frac{2E_T}{\MM_T}\right]
&\le\E\left[\left(f_Y(E_T)\MM_T^{1/q}\log\frac{2E_T}{\MM_T}\right)^p\right]^{1/p}\E\left[\MM_T^{q-1}\right]^{1/q}.
\end{align}
Note that for any $t>0$, the function
$h\mapsto h^{1/q}\log(2t/h)$ on $(0,\infty)$
attains its global maximum
$\frac{q}{e}(2t)^{1/q}$.
Therefore,
\begin{align}
\E\left[\left(f_Y(E_T)\MM_T^{1/q}\log\frac{2E_T}{\MM_T}\right)^p\right]^{1/p}
&\le
\frac{q\,2^{1/q}}{e}\E\left[f_Y(E_T)^p E_T^{p/q}\right]^{1/p}\notag\\
&\le C_{\varepsilon,K,\phi,{\gamma}}
(1+|x_0|^{2+\gamma})e^{C_{\varepsilon,K,\phi,{\gamma}}T},\label{1101-a}
\end{align}
where the last inequality follows from the inequalities in \eqref{1506+}
and \eqref{A.3}. 

On the other hand, since $q-1\in \mathbb{N}$, it follows from
inequality \eqref{1507+} with $q$ replaced by $q-1$ that
\begin{align}
\label{1101-b}
\E\left[\MM_T^{q-1}\right]^{1/q}\le\left(\sum_{i=0}^{L-1}\E[(\delta E_i)^{q-1}]\right)^{1/q}
\le C_q 
U(T)^{1/q}U(\Delta)^{1-2/q}.
\end{align}
Recall from \eqref{U-exp-bnd} and \eqref{U-and-phi} that
$U(T)\le C_\phi e^T$ and $U(\Delta)\le e/\phi(1/\Delta)$, respectively. Since
$2/q<\varepsilon$ and $\phi(1)/\phi(1/\Delta)\le1$ for
$\Delta\in(0,1)$,
\begin{align}\label{1101-c}
U(T)^{1/q}U(\Delta)^{1-2/q}\le C_{\varepsilon,\phi}e^T\left(\frac{\phi(1)}{\phi(1/\Delta)}\right)^{1-\varepsilon}.
\end{align}
Hence, combining \eqref{1506}-\eqref{1101-c} gives
\begin{align}\label{1113-b}
\E[I_2(T)]\le C_{\varepsilon,K,\phi,{\gamma}}(1+|x_0|^{2+\gamma})e^{C_{\varepsilon,K,\phi,{\gamma}}T}\left(\frac{1}{\phi(1/\Delta)}\right)^{1-\varepsilon}.
\end{align}
Finally, putting \eqref{1113-a} and \eqref{1113-b} together gives the desired result.
}
\end{proof}

\begin{remark}\label{r:rate1}
A remark similar to Remark \ref{r:rate0}(b) concerning Theorem \ref{thm:strong-approx} also applies to Theorem \ref{thm:synchro-strong}. Namely, 
the constants $C_q$ and $C_{\rho,q}$ in \eqref{1602-a} and \eqref{1602-b} diverge to $\infty$ as $\varepsilon\downarrow 0$, and consequently, the above proofs do not allow the $L^2$ convergence rate in \eqref{1607-new} to improve to 
$\phi(1/\Delta)^{-\min\{1,2\rho\}}$.
In fact, Section \ref{s:sharp} shows that such an improvement is generally impossible.
\end{remark}

\section{Remarks}
\label{s:remarks}

This section provides several noteworthy remarks about the statements derived in Sections \ref{s:approx-1} and \ref{s:approx-2}. In particular, Remarks \ref{r:compare1}--\ref{r:compare-deterministic-random} compare the convergence rates of our theorems with each other and with the existing literature. 
Remark \ref{r:square-int} discusses how the initial conditions of the time-changed SDEs can be randomized. 
The final remark, Remark \ref{r:t_0>0}, clarifies why our approach in Section \ref{s:approx-2} would not allow the initial time in SDE \eqref{SDE_002} to be positive.

\begin{remark}[{Comparing the rates of Theorems \ref{thm:strong-approx} and \ref{thm:synchro-strong}}]
\label{r:compare1}
Let the initial time be $t_0=0$, and let the SDE coefficients $F$ and $G$ be homogeneous; i.e., they depend on the space variable only. Then we can take $\rho=1$ in \eqref{SDE_condition4-a}, and both Theorems \ref{thm:strong-approx} and \ref{thm:synchro-strong} are applicable, provided conditions \eqref{SDE_condition1-a}--\eqref{SDE_condition2-a} are satisfied.

Theorem \ref{thm:strong-approx}, together with the inequality $\phi(\lambda)\le c \lambda$ valid for all $\lambda\ge 1$, implies that the rate of convergence is $\phi(1/\Delta)^{-1+\varepsilon}$ with any $\varepsilon\in(0,1)$. This coincides with the rate obtained by Theorem \ref{thm:synchro-strong}. 
In other words, although the approaches taken to proving Theorems \ref{thm:strong-approx} and \ref{thm:synchro-strong} are significantly different, the resulting rates for this sub-class of SDEs coincide.
\end{remark}

\begin{remark}[{Comparing the rates of Theorem \ref{thm:approx-synchr-conti} and \cite[Theorem 5.3]{BCT2026}}]
\label{r:compare2}
    In \cite{BCT2026}, a different approximation method of the time-changed It\^o diffusion was considered. 
    For the classical It\^o diffusion $Y$ in \eqref{SDE_002-Y},
    \cite[Theorem 5.3]{BCT2026} considers the Euler--Maruyama approximation $\widehat{Y}$ over the following time points including random times generated by the inverse subordinator $E$:
    \[
    \{0,\Delta,2\Delta,3\Delta,\ldots\}\cup \{E_{t_1},E_{t_2},\ldots\}.
    \]
Then an approximation process $\widehat X$ of the solution $X$ to SDE \eqref{SDE_002} is constructed with continuous interpolation so that it satisfies
$\widehat{X}_{t_i}=\widehat{Y}_{E_{t_i}}$, $i\ge 0$. Then a minor modification of the proof of \cite[Theorem 5.3]{BCT2026} yields
    \begin{align*}
         \E\left[\sup_{0\le s\le T}|X_s-\widehat X_s|^2\right]\le \widehat C e^{\widehat C T}\Delta^{\min\{1,2\rho\}},\qquad \Delta\in(0,1),
    \end{align*}
where $\widehat C=\widehat C(\phi,\rho,K, x_0)>0$.
Thus, this approximation method converges with order $\min\{\frac 12,\rho\}$, which is higher than the convergence order obtained from \eqref{1517} in Theorem \ref{thm:approx-synchr-conti}. 
    However, this approximation method tends to be resource-demanding. Indeed, for a fixed time horizon $T$, instead of considering the Euler--Maruyama approximation over the $L$ time points $E_{t_1},E_{t_2},\ldots,E_{t_L}=E_T$, one has to compute the approximation over $L+\lfloor\frac{E_T}{\Delta}\rfloor$ time points, which leads to two clear consequences. First, the fact that the number of time points is random may make implementation difficult. Second, taking a smaller step size $\Delta$ directly increases the number of time points. In other words, the higher rate in \cite[Theorem 5.3]{BCT2026} comes at the cost of a randomly determined, large number of discretization points.
\end{remark}

\begin{remark}[{Comparing the rates under deterministic time steps and those under random time steps}] \label{r:compare-deterministic-random}
Papers \cite{ChenLiuWu2026,DengLiu2020,Shen+2023}, which build upon the approaches taken in \cite{JinKobayashi2019,JinKobayashi2021,JumKobayashi},  
use approximation schemes that first simulate the underlying subordinator $D$ at equidistant deterministic time points $0=t_0<t_1<t_2<\ldots$ with $\Delta=t_{i+1}-t_i$ and then discretize the inverse subordinator $E$ at the random time steps
\begin{align}\label{1513}
    \tau_i\coloneqq D_{t_i}, \quad i\ge 0.
\end{align}
Although the use of the random time steps to approximate $E$ (and hence, approximate the solution of a time-changed SDE as well) may seem complicated, the analysis actually becomes cleaner with this approach since the increments of $E$ over the random time steps are deterministic; in fact,
\begin{align}\label{1514}
    E_{\tau_{i+1}}-E_{\tau_i}
    =t_{i+1}-t_i
    =\Delta, \quad i\ge 0.
\end{align}
In particular, this approach does not require statements on incremental moments of the inverse subordinator $E$, as in \eqref{Eincrements-stable} or Theorem \ref{Theorem_increments}. 

Suppose now that the SDE coefficients $F$ and $G$ satisfy conditions \eqref{SDE_condition1-a}--\eqref{SDE_condition4-a}. In \cite{JumKobayashi}, the solution $X$ to SDE \eqref{SDE_002} with synchronous coefficients is approximated by an Euler--Maruyama-type scheme over the random time steps $\tau_i$ in \eqref{1513}, with continuous interpolation. The proof of \cite[Theorem 3.2]{JumKobayashi} shows that the order of convergence for the uniform $L^1$ error in this setting is $\min(\frac{1-\varepsilon}{2},\rho)$ with any $\varepsilon\in(0,1)$. On the other hand, Theorem \ref{thm:approx-synchr-conti}, which approximates the solution $X$ over deterministic time steps in \eqref{deterministic-time-a}, shows that if $\phi$ satisfies a weak lower scaling condition with exponent $\alpha\in(0,1]$, then the corresponding order is $\frac{\alpha}{2}(\min\{1,2\rho\}-\varepsilon)$. The difference between the two orders becomes clearer if $\rho$ in \eqref{SDE_condition4-a} satisfies $\rho\ge \frac 12$, in which case, \cite[Theorem 3.2]{JumKobayashi} gives $\frac{1-\varepsilon}{2}$, whereas Theorem \ref{thm:approx-synchr-conti} gives $\frac{\alpha(1-\varepsilon)}{2}$. Consequently, we can conclude that the approximation based on the random time steps in \eqref{1513} leads to a higher convergence order than the approximation based on the deterministic time steps in \eqref{deterministic-time-a}. 
This difference is ascribed to the simple expression in \eqref{1514} that is available only when the random time steps in \eqref{1513} are employed. 

In contrast, \cite{JinKobayashi2019} analyzes a Euler--Maruyama-type approximation of SDE \eqref{SDE_001-a} with asynchronous coefficients, with piecewise-constant interpolation over the random time steps $\tau_i$ in \eqref{1513}. In this case, \cite[Theorem 7]{JinKobayashi2019} shows that the order of uniform $L^1$ convergence can be expressed as $\min\{\rho,\frac{1-\varepsilon}{2}\}$. Again, this is a higher order than $\min\{\rho,\frac{\alpha(1-\varepsilon)}{2}\}$ that follows from Theorem \ref{thm:strong-approx} when deterministic time steps in \eqref{deterministic-time-a} are used. 

However, let us point out that the validity of \cite[Theorem 7]{JinKobayashi2019} is limited to the case when (i) the Laplace exponent $\phi$ is regularly varying with index $\alpha\in(\frac 12,1)$, or (ii) $G(t,x)=G(t)$. This is because when the coefficient $G$ depends on $x$, the obtained upper bound for the uniform error contains a factor of the form $\E[e^{f(E_T)}]$ with $f$ being a polynomial of degree 2, and condition (i) is sufficient (but not necessary) to guarantee the finiteness of this expectation (recall Theorem \ref{theorem_expo_moment} (a)). 
Note that Theorem \ref{thm:strong-approx} does not have such limitations, as its proof used the method of exponential tilting to avoid the emergence of an expectation of the form $\E[e^{f(E_T)}]$. Utilizing this method in the proof of \cite[Theorem 7]{JinKobayashi2019} allows the limitations of the theorem to be removed.
\end{remark}

\begin{remark}[{Discussion of random initial values}]\label{r:square-int}
In Sections \ref{s:approx-1} and \ref{s:approx-2}, we took the initial values of the SDE solutions to be deterministic so as to highlight the core ideas behind the proofs. However, this restriction can be
relaxed when the initial values are independent of the driving processes.
\begin{enumerate}[(a)]
\item In Theorem \ref{thm:strong-approx}, we may replace the initial
value $x_{t_0}$ with a square-integrable random variable $X_{t_0}$ that is independent of $(B,E)$. 
To see this, recall that the random variable $H_{q,T}$ in the proof of Theorem \ref{thm:strong-approx} 
 involves $\E_B[Y_T^{(2)}]$ and $\left(\E_B[Y_T^{(2q)}]\right)^{1/q}$. It follows from the conditional version of Lemma \ref{lemma_Holder-a} together with the simple inequality $(1+|x_{t_0}|^{2q})^{1/q}\le1+|x_{t_0}|^2$ that
\[
\E_B[Y_T^{(2)}]\le C_K(1+|x_{t_0}|^2)e^{C_KE_T},\qquad\left(\E_B[Y_T^{(2q)}]\right)^{1/q}\le C_{q,K}(1+|x_{t_0}|^2)e^{C_{q,K}E_T}.
\]
Consequently, the inequality in \eqref{HqT} can be re-expressed as
\[\E[H_{q,T}]\le \left(\E[H_{q,T}^p]\right)^{1/p}\le C_{q,K,\phi}(1+|x_{t_0}|^2)e^{C_{q,K,\phi}T},
\]
where $p=q/(q-1)$ and the constant $C_{q,K,\phi}$ is independent of $x_{t_0}$. With the remaining estimates in the proof unchanged, the inequality in \eqref{1427} holds with the initial-value factor $1+|x_{t_0}|^2$.

As a result, replacing $x_{t_0}$ with the above-mentioned random variable $X_{t_0}$ and using the conditioning argument gives
\begin{align*}
\E\left[\sup_{t_0\le t\le T}|X_t-X_t^\Delta|^2\right]\le C(1+\E[|X_{t_0}|^2])e^{CT}\left(\Delta^{2\rho}+\phi(1/\Delta)^{-(1-\varepsilon)}\right),\quad \Delta\in(0,1),
\end{align*}
where $C=C(K,\phi,\varepsilon)$ is independent of $T$, $\Delta$, the partition, and the distribution of $X_{t_0}$. 

\item Theorem \ref{thm:synchro-strong} may also be written with a square-integrable random initial value $X_0$ that is independent of $(B,E)$. However, a careful examination of the proof of Theorem \ref{thm:synchro-strong} shows that the bound in \eqref{1607-new} can be written as
\[
    \E\left[\sup_{0\le s\le T}|X_s- X^\Delta_s|^2\right]\le 
        {C(1+|x_0|^2)(1+|x_0|^{2+\gamma})e^{CT}\phi(1/\Delta)^{-(\min\{1,2\rho\}-\varepsilon)}},
\]
for some $C=C(\phi,\varepsilon,K,\rho)>0$ and $\gamma>0$. In other words, the bound contains the $(2+\gamma)$-th power of $x_0$ which comes from the bound for $I_2(T)$ defined in \eqref{I1+I2}. We therefore modify the estimate of $I_2(T)$ in \eqref{I1+I2} to obtain a bound depending only on the factor $1+|x_0|^2$.

To this end, fix $\varepsilon\in(0,\min\{1,2\rho\})$, choose an integer $q>1/\varepsilon$, and put $p=q/(q-1)$. As before, in what follows, constants  may vary from line to line. For $0\le a<a+h\le S$, the Burkholder--Davis--Gundy inequality in \eqref{Burkholder}, H\"older's inequality, the linear-growth bound \eqref{SDE_condition2-a}, and the supremum bound \eqref{1314a} with exponent $2q$ give
\begin{align*}
\E_B\left[\sup_{a\le r\le a+h}|Y_r-Y_a|^{2q}\right]
&\le C_{q,K}(h^{2q-1}+h^{q-1})
\int_a^{a+h}\E_B[1+|Y_r|^{2q}]\,\dd r\\
&\le C_{q,K}(1+|x_0|^{2q})e^{C_{q,K}S}h^q,
\end{align*}
where $h\le S$ allows the factor $1+h^q$ to be absorbed into the exponential. Note that the bound also holds for $h=0$.
Since $B$ and $E$ are independent, the preceding bound applies with $S=E_T$, $a=E_{t_i}$, and $h=\delta E_i$. Using $(1+|x_0|^{2q})^{1/q}\le1+|x_0|^2$ and the definition of $I_2(T)$ in \eqref{I1+I2}, we obtain
\begin{align*}
\E_B[I_2(T)]
&\le 2\left(\sum_{i=0}^{L-1}\E_B\left[
\sup_{E_{t_i}\le r\le E_{t_{i+1}}}|Y_r-Y_{E_{t_i}}|^{2q}
\right]\right)^{1/q}\\
&\le C_{q,K}(1+|x_0|^2)e^{C_{q,K}E_T}
\left(\sum_{i=0}^{L-1}(\delta E_i)^q\right)^{1/q}.
\end{align*}
Taking the expectation with respect to $E$ and applying H\"older's inequality together with the bounds in  \eqref{A.3} and \eqref{1507+} gives
\begin{align}
\E[I_2(T)]
&\le C_{q,K}(1+|x_0|^2)
\E[e^{pC_{q,K}E_T}]^{1/p}
\left(\sum_{i=0}^{L-1}\E[(\delta E_i)^q]\right)^{1/q}\notag\\
&\le C_{q,K,\phi}(1+|x_0|^2)e^{C_{q,K,\phi}T}
U(T)^{1/q}U(\Delta)^{1-1/q}\notag\\
&\le C_{\varepsilon,K,\phi}(1+|x_0|^2)e^{C_{\varepsilon,K,\phi}T}
\phi(1/\Delta)^{-1+\varepsilon}.\label{1135-newbound}
\end{align}
For the last inequality, we used $U(\Delta)\le U(T)$ and $1/q<\varepsilon$ to bound the product of powers of $U$ by $U(T)^\varepsilon U(\Delta)^{1-\varepsilon}$, and then apply \eqref{U-exp-bnd} and \eqref{U-and-phi}.
Combining the estimate in \eqref{1135-newbound} with the estimate of $I_1$ in \eqref{1113-a} gives \eqref{1607-new} with the initial-value factor depending only on $1+|x_0|^2$.

As a result, replacing $x_{t_0}$ with a random variable $X_{t_0}$ independent of $(B,E)$ gives
\begin{align*}
\E\left[\sup_{0\le s\le T}|X_s-X_s^\Delta|^2\right]
\le C(1+\E[|X_0|^2])e^{CT}
\phi(1/\Delta)^{-(\min\{1,2\rho\}-\varepsilon)},
\quad \Delta\in(0,1),
\end{align*}
where $C=C(\phi,\varepsilon,K,\rho)$ is independent of $T$, $\Delta$, the partition, and the distribution of $X_0$.
\end{enumerate}

\end{remark}

\begin{remark}[{Discussion of why the initial time in SDE \eqref{SDE_002} must be 0}]
\label{r:t_0>0}
    Fix a time $t_0>0$, and let  \(T\coloneqq E_{t_0}\).  Then $T>0$ a.s. Consider the SDE
\begin{align*}
		X_t
		&=x_{t_0}+\int_{t_0}^{t} F(E_s,X_s)\,\dd E_s+\int_{t_0}^{t} G(E_s,X_s)\,\dd B_{E_s}\\
		&=x_{t_0}+\int_{T}^{E_t} F(s,X_{D_s})\,\dd s+\int_{T}^{E_t} G(s,X_{D_s})\,\dd B_s, \ \ t\ge t_0,
\end{align*}
where we used the change of variables formula in \cite[Lemma 2.3]{Kobayashi2011}. The filtration considered here is generated by $B$ and $D$ as follows:
\[
	\mathcal{F}_t\coloneqq \sigma(B_s, D_s;\, s\le t),	
\]
which is completed and made right-continuous so that it satisfies the usual conditions. Then the random time $T=E_{t_0}$ is an $(\mathcal{F}_t)$-stopping time, and the process $Y_\tau\coloneqq X_{D_\tau}$ satisfies 
\[
	Y_\tau=x_{t_0}+\int_{T}^{\tau} F(s,Y_s)\,\dd s+\int_{T}^{\tau} G(s,Y_s)\,\dd B_s, \ \ \tau\ge T.
\]
Define the following time-shifted objects:
\begin{align*}
\widetilde{B}_\tau &= B_{T+\tau}-B_{T}, &
\widetilde{D}_\tau &= D_{T+\tau}-D_{T}, &
\widetilde{\mathcal{F}}_\tau &= \mathcal{F}_{T+\tau}, \\
\widetilde{Y}_\tau &= Y_{T+\tau}, &
\widetilde{F}(\tau,z) &= F(T+\tau,z), &
\widetilde{G}(\tau,z) &= G(T+\tau,z).
\end{align*}
Since $B$ and $D$ are independent, the pair $(B,D)$ is a L\'evy process satisfying the strong Markov property with respect to $(\mathcal{F}_t)$. Therefore, $(\widetilde{B},\widetilde{D})$ is independent of $\mathcal{F}_T$ and is a L\'evy process adapted to $(\widetilde{\mathcal{F}}_\tau)$ and having the same law as $(B,D)$. In particular, $\widetilde{B}$ is a Brownian motion, $\widetilde{D}$ is a subordinator, and they are independent. 
It follows that 
\begin{align*}
	\widetilde{Y}_\tau
	&=x_{t_0}+\int_{T}^{T+\tau} F(s,Y_s)\,\dd s+\int_{T}^{T+\tau} G(s,Y_s)\,\dd B_s\\
	&=x_{t_0}+\int_{0}^{\tau} F(T+s,Y_{T+s})\,\dd s+\int_{0}^{\tau} G(T+s,Y_{T+s})\,d\widetilde{B}_s\\
	&=x_{t_0}+\int_{0}^{\tau} \widetilde{F}(s,\widetilde{Y}_s)\,\dd s+\int_{0}^{\tau} \widetilde{G}(s,\widetilde{Y}_s)\,d\widetilde{B}_s, \ \ \tau\ge 0.
\end{align*}
Now, define the inverse of the shifted subordinator $\widetilde{D}$ by 
\[
	\widetilde{E}_\tau\coloneqq \inf\{u>0: \widetilde{D}_u>\tau\}, \ \ \tau\ge 0. 
\]
Then the duality principle in \cite{Kobayashi2011} connects the SDE for $\widetilde{Y}$ with the SDE 
\[
	\widetilde{X}_\tau=x_{t_0}+\int_{0}^{\tau} \widetilde{F}(\widetilde{E}_s,\widetilde{X}_s)\,d\widetilde{E}_s+\int_{0}^{\tau} \widetilde{G}(\widetilde{E}_s,\widetilde{X}_s)\,d\widetilde{B}_{\widetilde{E}_s}, \ \ \tau\ge 0,
\]
via the relationships $\widetilde{X}_\tau=\widetilde{Y}_{\widetilde{E}_\tau}$ and $\widetilde{Y}_\tau=\widetilde{X}_{\widetilde{D}_\tau}$. The latter SDE with initial time 0 is intrinsically different from the time-changed SDE in \eqref{SDE_002}. In fact, the coefficient functions are random since $\widetilde{F}(u,x)=F(T+u,x)=F(E_{t_0}+u,x)$ and $\widetilde{G}(u,x)=G(T+u,x)=G(E_{t_0}+u,x)$. Furthermore, 
\begin{align*}
	\widetilde{E}_\tau
	&=\inf\{u: u>0, D_{T+u}-D_T>\tau\}
	=\inf\{r: r>0, D_r>D_T+\tau\}-T=E_{D_T+\tau}-T,
\end{align*}
which yields
\[
	\widetilde{X}_\tau
=\widetilde{Y}_{\widetilde{E}_\tau}=Y_{T+\widetilde{E}_\tau}=Y_{E_{D_T+\tau}}
	=X_{D_T+\tau}
	=X_{D_{E_{t_0}}+\tau}
    \ne X_\tau.
\]
This discussion demonstrates that the case when the initial time $t_0$ is positive cannot be handled within the framework of Section \ref{s:approx-2}. 
\end{remark}

\section{Sharpness of the convergence rate}\label{s:sharp}

Remarks \ref{r:rate0} and \ref{r:rate1} observed that the constants appearing in the error bounds of Theorems \ref{thm:strong-approx} and \ref{thm:synchro-strong} diverge to $\infty$ as $\varepsilon\downarrow0$, and consequently, the methods employed to prove Theorems \ref{thm:strong-approx} and \ref{thm:synchro-strong} do not allow the $L^2$ convergence rate to improve to $\Delta^{2\rho}+\phi(1/\Delta)^{-1}$ and $\phi(1/\Delta)^{-\min\{1,2\rho\}}$, respectively. 
Both of these expressions reduce to the form $c\Delta^\alpha$ if the H\"older exponent in \eqref{SDE_condition4-a} is $\rho=1$ and the Laplace exponent $\phi$ satisfies $\phi(\lambda)\ge c\lambda^\alpha$ for all $\lambda\ge\lambda_0$.
However, the latter discussion itself does not conclude that achieving $\Delta^\beta$ with $\beta\ge \alpha$ as a rate function for the $L^2$ strong error is impossible.  
The following theorem, in contrast, confirms that such an improvement is generally impossible when piecewise-constant interpolation is used for approximation, thereby confirming the sharpness of the results presented in Theorems \ref{thm:strong-approx} and \ref{thm:synchro-strong}.

\begin{theorem}\label{prop:sharp-piecewise}
Let $E$ be the inverse of an $\alpha$-stable subordinator $D$ with $\alpha\in(0,1)$. Let $B$ be a Brownian motion independent of $E$. 
Consider SDEs \eqref{SDE_001-a} and \eqref{SDE_002} with coefficients $F\equiv 0$ and $G\equiv 1$, so that the SDE solution and its approximation with piecewise-constant interpolation are respectively given by $X_t= B_{E_t}$ and $X_t^\Delta=B_{E_{t_{\ell_t}}}$ for $t\in[0,T]$, where the equidistant mesh with $\Delta=T/L$ is taken; i.e., $t_\ell=\ell \Delta$ for all $\ell$ in \eqref{deterministic-time-a}. Then
\[
\lim_{\Delta\downarrow0}
\frac{\E[\sup_{0\le t\le T}|X_t-X_t^\Delta|^2]}{\Delta^\alpha}
=\infty.
\]
\end{theorem}

The proof of Theorem \ref{prop:sharp-piecewise} requires the following simple lemma, which, to our knowledge, has not appeared in the literature. 

\begin{lemma}\label{lemma:max_increments}
Let $E$ be the inverse of a subordinator $D$. 
Suppose that there exists a positive sequence $\{a_n\}$ such that $a_n\to\infty$ as $n\to\infty$ and $\P(D_{a_n}\le 1)>0$ for all $n$. Let
\[
G_L\coloneqq \max_{0\le i\le L-1} (E_{i+1}-E_i), \quad L=1,2,3,\ldots. 
\]
Then a.s., $G_L \to \infty$ as $L\to \infty$.

\end{lemma}

    \begin{proof}
    For a fixed \(n\in\N\), define
    \[
        A_k^{(n)}
        \coloneqq 
        \left\{
            D_{(k+1)a_n}-D_{ka_n}\leq1
        \right\}, \quad k\in\N\cup \{0\}.
    \]
    Since \(D\) has stationary and independent increments, \(\{A_k^{(n)}\}_{k\geq0}\) are independent, and 
    \(
        \mathbb{P}\bigl(A_k^{(n)}\bigr)
        =
        \mathbb{P}(D_{a_n}\leq1)
        >0
    \)
    for all $k$.
    Therefore, the second Borel--Cantelli lemma yields
    \begin{equation}
        \mathbb{P}\bigl(A_k^{(n)}\text{ occurs for infinitely many $k$}\bigr)=1.
        \label{3}
    \end{equation}
    
    Now, on the event \(A_k^{(n)}\),
\(
    D_{(k+1)a_n}\leq D_{ka_n}+1\leq m_k+2,
\)
where $m_k\coloneqq \lceil D_{ka_n}\rceil-1$.
Since $D_t\to\infty$ as $t\to\infty$ a.s., we restrict our attention to sufficiently large $k$ so that $m_k\ge1$.
    The latter implies that the range of the path of 
    \(
        \bigl(D_r\bigr)_{r\in[ka_n,(k+1)a_n]}
    \)
    is contained in the interval \([m_k,m_k+2)\). Consequently,
    \[
        [ka_n,(k+1)a_n]
        \subset
        [E_{m_k},E_{m_k+2}].
    \]
    Hence, 
    \begin{equation}\label{4}
        a_n
        \leq
        E_{m_k+2}-E_{m_k}
        =
        \sum_{j=1}^{2}
        \bigl(E_{m_k+j}-E_{m_k+j-1}\bigr)\le 2\max_{1\leq j\leq2}
        \bigl(E_{m_k+j}-E_{m_k+j-1}\bigr).  
    \end{equation}
    The inequality in \eqref{4} implies that there exists at least one $i_k^{(n)}\in\N$ such that $E_{i_k^{(n)}+1}-E_{i_k^{(n)}}\ge \frac{a_n}{2}$, and hence, $G_L\ge \frac{a_n}{2}$ for all $L\ge i_k^{(n)}+1$. Therefore, $\liminf_{L\to\infty}G_L\ge \frac{a_n}{2}$. 
    Thus, for each fixed $n\in\N$, by \eqref{3}, $\liminf_{L\to\infty}G_L\ge \frac{a_n}{2}$ a.s. The latter implies that, a.s., $\liminf_{L\to\infty}G_L\ge \frac{a_n}{2}$ for all $n\in\N$, which yields the desired result.
    \end{proof}

\begin{proof}[Proof of Theorem \ref{prop:sharp-piecewise}]
Let
\[
    \delta E_i\coloneqq E_{t_{i+1}}-E_{t_{i}}=E_{(i+1)\Delta}-E_{i\Delta},
    \qquad i=0,\ldots,L-1.
\]
Let
\begin{align*}
    \mathcal{W}_{\Delta}(T)
    &\coloneqq 
    \sup_{0\leq t\leq T}
|X_t-X_t^\Delta|
    =\sup_{0\leq t\leq T}
    \left|B_{E_t}-B_{E_{t_{\ell_t}}}\right|\\
    &=\max_{0\leq i\leq L-1}
    \sup_{t_i\leq t\leq t_{i+1}}
    \left|B_{E_t}-B_{E_{t_i}}\right|\\
    &=
    \max_{0\leq i\leq L-1}
    \sup_{0\leq r\leq \delta E_i}
    \left|B_{E_{t_{i}}+r}-B_{E_{t_{i}}}\right|.
\end{align*}
Conditionally on \(E\), the increments of the Brownian motion over the disjoint
intervals $[E_{t_{i}},E_{t_{i+1}})$, $i=0,\ldots,L-1$, are independent. 
On the other hand, the Markov property of the Brownian motion implies that, given $E$, 
\[
    (B_{E_{t_{i}}+r}-B_{E_{t_{i}}})_{r\in[0, \delta E_i]}
     \stackrel{d}{=}
     (B_r)_{r\in[0, \delta E_i]},
\]
and hence, by Brownian scaling,
\[
    \sup_{0\leq r\leq \delta E_i}
    \left|B_{E_{t_{i}}+r}-B_{E_{t_{i}}}\right|
    \stackrel{d}{=}
    \sup_{0\leq r\leq \delta E_i}|B_r|    
    \stackrel{d}{=}
    (\delta E_i)^{1/2}
    \sup_{0\leq r\leq 1}|B_r|.
\]
Now, take independent copies $B^{(0)},\ldots,B^{(L-1)}$ of the Brownian motion $B$ that are independent of the inverse subordinator $E$, and let 
\[
    S_i\coloneqq \sup_{0\leq r\leq 1}|B^{(i)}_r|,
    \qquad i=0,\ldots,L-1.
\]
Then the above discussion implies that, conditionally on $E$,
\begin{align}\label{eq:sharp:goal2}
    \mathcal{W}_{\Delta}(T)
    \stackrel{d}{=}
    \max_{0\leq i\leq L-1}(\delta E_i)^{1/2}S_i.
\end{align}

Next, for each $L\ge 1$, 
let
\[
    I_L
    \coloneqq 
    \min\left\{
        \ell\in\{0,\ldots,L-1\}:
        \delta E_\ell=\max_{0\le i\le L-1} \delta E_i
    \right\}
\]
be the least index $\ell$ for which the maximum is attained.
Then 
\begin{align}\label{eq:sharp:goal3}
\max_{0\le i\le L-1} (\delta E_i)S_i^2
\ge (\delta E_{I_L})S_{I_L}^2
= S_{I_L}^2 \max_{0\le i\le L-1} \delta E_i.
\end{align}
By the conditional distributional equality and inequality in \eqref{eq:sharp:goal2} and \eqref{eq:sharp:goal3},
\begin{align*}
    \mathbb{E}\bigl[(\mathcal{W}_\Delta(T))^2 | E\bigr]
    \geq
    \E[S_{I_L}^2|E]\max_{0\le i\le L-1} \delta E_i \quad \textrm{a.s.}
\end{align*}
Note that $\E[S_{I_L}^2|E] = \E\left[\sup_{0\le t\le 1}B_t^2\right]$ a.s., which is finite due to Doob's maximal inequality.
Thus, taking the expectation with respect to $E$ gives 
\begin{align}\label{eq:sharp:goal5}
    \mathbb{E}\bigl[(\mathcal{W}_\Delta(T))^2\bigr]
    \geq
    \E\left[\sup_{0\le t\le 1}B_t^2\right]\E\left[\max_{0\le i\le L-1} \delta E_i\right].
\end{align}
Since the inverse $\alpha$-stable subordinator is self-similar with index $\alpha$,
$(E_{t_i})_{0\le i\le L}\stackrel{d}{=}(\Delta^\alpha E_i)_{0\le i\le L}$, which implies that 
$(\delta E_i)_{0\le i\le L-1}\stackrel{d}{=}\left(\Delta^\alpha (E_{i+1}-E_i)\right)_{0\le i\le L-1}$. Therefore, 
\begin{align}\label{eq:sharp:goal6}
    \E\left[\max_{0\le i\le L-1} \delta E_i\right]
    =\Delta^\alpha \E[G_L],
\end{align}
where $G_L$ is defined in Lemma \ref{lemma:max_increments}.

Finally, combining \eqref{eq:sharp:goal5} and \eqref{eq:sharp:goal6} gives 
\begin{equation}
    \frac{\mathbb{E}\bigl[(\mathcal{W}_\Delta(T))^2\bigr]}{\Delta^\alpha}
    \geq
    \E\left[\sup_{0\le t\le 1}B_t^2\right]\mathbb{E}[G_L].
    \label{7}
\end{equation}
Notice that the stable subordinator $D$ satisfies the assumption of Lemma \ref{lemma:max_increments} since \(D_n\) for any $n\in\N$ has a strictly positive density on \((0,\infty)\); see, e.g., \cite[Remark 4 in §2.2]{Zolotarev}.
Moreover, $G_L$ is increasing with respect to $L$, and $L\to\infty$ as $\Delta\downarrow 0$. Hence, the inequality in \eqref{7}, the monotone convergence theorem,  and Lemma \ref{lemma:max_increments} together yield the desired result. 
\end{proof}

\begin{remark}
Theorem \ref{prop:sharp-piecewise} shows that the convergence orders obtained in Theorems \ref{thm:strong-approx} and \ref{thm:synchro-strong} cannot be improved to the case $\eps=0$ in general. However, the identification of the exact rate is still an open question. The exact rate may involve a logarithmic correction, and the following two known facts from the literature are worth mentioning in this regard. {First, \cite[Equation (10)]{FischerNappo} presents a logarithmic correction factor under the untime-changed setting; see also \cite[Lemmas 1 and 2]{FischerNappo}.} Second, \cite[Lemma 3.4] {OS15} shows that the random variable $M$ defined by  
    \begin{align}\label{eq:1309}
        M\coloneqq \sup_{0 \leq v < u \leq 1/2}
\frac{E_u-E_v}
{(u-v)^\alpha \lvert \log(u-v)\rvert^{1-\alpha}}
    \end{align}
    has finite moments of all orders.
    The estimate \eqref{eq:1309} seems to suggest that the maximal increment $\max_{0\leq i\leq L-1}\delta E_i$ of the inverse stable subordinator is, up to a random factor with finite moments, bounded above by a quantity of order $\Delta^\alpha \lvert\log \Delta\rvert^{1-\alpha}$. Since $\mathcal{W}_{\Delta}(T)$ appearing in the proof of Theorem \ref{prop:sharp-piecewise} is the modulus of continuity of Brownian motion evaluated over these random time increments, a logarithmic correction may be needed. 
\end{remark}

\section*{Acknowledgments}
Ivan Biočić acknowledges financial support by the European Union – NextGenerationEU through the National Recovery and Resilience Plan 2021-2026 Institutional grant of University of Zagreb Faculty of Science (IK IA 1.1.3. Impact4Math), as well as the support by Croatian Science Foundation through the project IP-2025-02-8793.

Part of this research was conducted while Kei Kobayashi was on a Faculty Fellowship at Fordham University. He appreciates support from Fordham University.

\newpage
\bibliographystyle{abbrv}
\bibliography{references}

\end{document}